\documentclass{amsart}
\usepackage[english]{babel}
\pdfoutput=1 

\usepackage{all2021}                    
\usepackage{tree-commands4}             

\usepackage{overpic}                    
\usepackage{tikz-cd}                    
\newcommand{\figsdir}{figs}

\usepackage[colorlinks, citecolor={black}, linkcolor={{black}}, urlcolor={blue}]{hyperref}          
\calclayout

\usepackage{enumerate}

\usepackage[
			backend=biber, 
			natbib=true,
			bibencoding=utf8,
			style=alphabetic,
			labelalpha=true,
			maxalphanames=5, 
			minalphanames=2,	
]{biblatex}  
\usepackage{csquotes}

\title[Combinatorics of higher-dimensional tropical covers]{Combinatorics of higher-dimensional tropical covers}
\author{Alejandro Vargas}

\begin{document} 

		\begin{abstract}
We develop a combinatorial framework to study certain polyhedral maps which are higher-dimensional analogues of tropical covers between metric graphs. 
Under a mild combinatorial assumption, we show that a map satisfies the so-called balancing condition if and only if it is an indexed branched cover, i.e.~locally over connected sets the count with multiplicity of points in every fibre is a constant, which in particular gives a well-defined global degree when the target is connected. 
Given a balanced map $(\Sigma, m_\Sigma) \to \Delta$, we lift several connectivity properties of $\Delta$ to~$\Sigma$. 
Using these lifting results we determine whether a multiplicity $m_{\mathcal U}$ that is defined only on the interiors of maximal cells of $\Sigma$ can be extended to all $\Sigma$ in a balanced manner. 
This relies on a strong connectivity assumption; we give a counterexample when this is missing.

	\end{abstract}

		\maketitle

\section{Introduction}

\subsection{Counting points in fibres} Let $F : X \to Y$ be a continuous map. 
Recall that $F$ is a \emph{cover} if every point $y \in Y$ has a neighbourhood $U$ whose preimage is a disjoint union of copies of~$U$.
See Figure \ref{fig:covers}~(a) for an example.
The cardinality of the fibre $\inv{F}(y)$ is the \emph{degree} of $F$ at $y \in Y$. 
When $Y$ is connected, the degree is  an invariant of~$F$ independent of~$y$.
As a slight generalization,
for us a branched cover is a map $F : X \to Y$ such that $F$ is a cover over a dense and open subset of $Y$.
The complement of that set is called the \emph{branch locus} $B$.
See Figure \ref{fig:covers}~(b).
To link the degree over $B$ with $\inv F(Y \setminus B) $ and get an invariant, 
consider an index map $m_F: X \to \ZZ_{\ge 1}$, and count points in $X$ with multiplicity~$m_F$.
We say that $(F, m_F)$ is an \emph{indexed branched cover} if for any connected open $U \subset Y$ and any connected component $V$ of $\inv F(U)$, the map  $F$ restricts to a map with a constant degree.
Thus, we can speak of local degrees, and if $Y$ is connected, once again we get a well-defined constant $\deg F$.
See Figure \ref{fig:covers}~(c) for an example, and Definitions~\ref{de:local-degree} and~\ref{de:IndexedBranchedCover} for details.

This definition is motivated by maps of Riemann surfaces, tropical morphisms of metric graphs, and maps of tropical moduli spaces. These objects have natural index maps, although in the case of tropical moduli spaces often this multiplicity is defined just on the dense subset of interiors of top-dimensional cells, which brings us to our first question.

\begin{que} 
    \label{que:Extension}
    Let $F: X \to Y$ be a branched cover, and $m_V$ a map defined on a dense subset $V$ of $X$.
    When does $m_V$ extend to a map $m_F$ such that $(F, m_F)$ is an indexed branched cover?
\end{que}

Given that the definition of indexed branched cover involves checking a condition over all connected open sets $U \subset Y$ and connected components $V$ of $\inv F(U)$, we desire tools that simplify this work.

\begin{que} 
    \label{que:Criteria}
    Establish criteria to check when a pair $(F, m_F)$ is an indexed branched cover. 
\end{que}

Finally, as with topological covers, the following question is relevant when investigating deformation arguments on $X$.

\begin{que} 
    \label{que:Connectivity}
    Let $(F: X \to Y, m_F)$ be an indexed branched cover. What connectivity properties of $Y$ lift to $X$?
\end{que}

\subsection{Results and structure of paper}
These questions are too broad for all topological spaces, 
so our efforts are directed by tropical geometry, 
i.e.~the study of piecewise-linear combinatorial shadows of phenomena in algebraic geometry.
Thus, we study abstract polyhedral complexes $\topo \Sigma$ and their face posets $\Sigma$. 
The latter is the quotient space obtained by identifying all points in the interior of a cell. 
So we have a natural projection $\poly_\Sigma: \topo \Sigma \to \Sigma$, 
which is continuous, surjective and also open. We require finite $\Sigma$.

These spaces are abstract because they arise by gluing polyhedra, so an embedding into a real vector space is not guaranteed. 
All the information is stored in a functor $\sigma: \Sigma \to \polycat$, 
and the topological realization $\topo \Sigma$ is obtained via a colimit.
We require the functor $\sigma$ to satisfy several lifting conditions to ensure that every polyhedral face and inclusion of $\topo \Sigma$ has a corresponding object or morphism in~$\Sigma$.
Our definition of polyhedral space follows closely \cite{acp15, cmr16, ccuw20}, and can handle self-gluing at faces and automorphisms. 
After making a few remarks about the general case, we assume that $\Sigma$ as a category is a poset, which gives abstract polyhedral complexes. 
Closed and open sets in $\Sigma$ are downward closed and upwards closed sets, generated by finite unions of principal down-sets $\downset \alpha$ and principal up-sets $\upset \alpha$, respectively.

Throughout this article we reduce questions about polyhedral complexes $\topo \Sigma$ to combinatorial questions about the face poset $\Sigma$; the main results of this nature being: 

    \begin{mytheorem} 
      \label{thm:SummarySection1}
     Let $\Phi : [\Sigma \to \polyintface] \to [\Delta \to \polyintface]$ be a morphism of polyhedral complexes, 
     $\varphi$ the induced morphism of posets, i.e.~$\dtmor$ is order preserving, 
     and $m_\varphi : \Sigma \to \ZZ_{\ge 1}$ an index map. 
     Assume that $\varphi$ maps down-sets isomorphically onto their image. 
     We have 
   \begin{enumerate}
         \item Diagram~\ref{dia:FiberProduct} is a fibre product in the category $\Top$ of topological spaces. 

      \begin{minipage}{0.9\linewidth}
			\centering
     \[\begin{tikzcd} 
        \topo \Sigma \arrow[r, "\poly_\Sigma"] \arrow[d, swap, "\topo \Phi"] & \Sigma \arrow[d, "\varphi"] \\
        \topo \Delta \arrow[r, "\poly_\Delta"] & \Delta.
      \end{tikzcd} \]
			\refstepcounter{thm}
			      \label{dia:FiberProduct}
			Diagram~\ref{dia:FiberProduct}
     \end{minipage}
     \vspace{1em}

     \item  If one of the pairs $(\dtmor, m_\dtmor)$ or $(\topo \Phi, \poly_\Sigma \circ m_\dtmor)$ is an indexed branched cover, then so is the other. 
     \end{enumerate}
    \end{mytheorem}

    The assumption on $\dtmor$ has a geometric origin and ensures that $\dtmor$ reflects properly the combinatorial information of $\topo \Phi$. 
    We say such morphism of posets is \combinatorial{}.    
    See Subsection~\ref{sub:CombinatorialMorphisms} and the discussion after Equation~\eqref{eq:containment-fibres} on Page~\pageref{eq:containment-fibres}.
    We argue in Subsection~\ref{sub:RefiningMorphisms} that this is a mild condition which can always be achieved with refinements.
    This result constitutes progress on Question~\ref{que:Criteria}, and its proof requires Lemmas~\ref{lm:LocalDegreeUnionTwoSets} and \ref{lm:PathConnectedCountablyBasisSimplifies}, themselves criteria for being an indexed branched cover. 
    In light of Theorem~\ref{thm:SummarySection1}, we now focus on finite posets.

    We generalize the well-known balancing condition of so-called tropical morphisms or tropical covers;
see e.g~\cite{mik07, bn09, bbm11, cap14} for background, or the expository~\cite{dv21}.
Given a morphism $\dtmor : \Sigma \to \Delta$ of posets
and an open set~$\calV \subset \Sigma$, 
we say that a map $m_\calV : \calV \to \ZZ_{\ge 1}$ is balanced with respect to $\dtmor$ if 
    for any $\alpha$ in $\calV$ and any choice of $\beta$ in $\Delta$ such that $\varphi(\alpha) \lessdot \beta$ we have that
    \begin{align} \label{eq:BalancingConditionIntro}
       m_\calV(\alpha) = \sum_{ \substack{ \gamma \in \inv \dtmor(\beta) \\ \alpha \lessdot \gamma} } m_\calV(\gamma).
    \end{align} 
    Here $\alpha \lessdot \gamma$ means that \emph{$\gamma$ covers $\alpha$}, i.e.~$\gamma$ is a minimum in the set $\calU_\alpha = \upset \alpha \setminus \aset{\alpha}$.
    In a polyhedral complex this means that $\sigma_\alpha \subset \sigma_\gamma$ and that $\dim \sigma_\alpha = \dim \sigma_\gamma - 1$.
	\begin{mytheorem} 
     \label{mytheorem:IndexedBranchedCovervsBalanced} 
		Let $\varphi : \Sigma \to \Delta$ be a \combinatorial{} morphism of posets,
    and $m_\varphi : \Sigma \to \ZZ_{\ge 1}$ an index map. 
    Assume that $\Delta$ is connected.
    The pair $(\varphi, m_\varphi)$ is an indexed branched cover with branch locus in the complement of  $\max \Delta$
    if and only if $(\varphi, m_\varphi)$ is balanced.
	\end{mytheorem} 
    Again, $\dtmor$ being \combinatorial{} is necessary: see Remark~\ref{re:CounterexamplesToIBCiffBalanced} for counterexamples in both directions when the condition is missing.
    Having the balancing condition as a tool, we explore Question~\ref{que:Connectivity} in Subsection~\ref{sub:ANecessaryConditionForQuestion}.
    The outcome is several results on lifting paths. 
    First, Lemma~\ref{lm:LiftingUpwardsPaths}  lifts upwards paths, which implies a necessary condition for Question~\ref{que:Extension}.

\begin{myprop} 
    \label{myprop:IndexedBranchedCoversAreOpenMaps}  
     Let $\varphi: \Sigma \to \Delta$ be a morphism of posets.
     If there exists a balanced map $m_\Sigma : \Sigma \to \ZZ_{\ge 1}$,
     then $\varphi$ is an open map. 
\end{myprop}

   On the side of a sufficient condition, we initially conjectured that it was enough to have that $\Sigma$ is connected through codimension-1.
   This is false, see Example~\ref{ex:CannotAlwaysExtend} that has a map $\dtmor: \Sigma \to \Delta$ where $\Sigma$ is connected through codimension-1, $\dtmor$ satisfies the necessary condition of being open, yet remarkably we exhibit an index map which cannot be extended. 
   We are able to prove extension under a stronger connectivity condition.

\begin{mytheorem} 
   \label{mytheorem:ExtendIndexMap}
  Let $\varphi : \Sigma \to \Delta$ be a \combinatorial{} morphism 
  and $\calV  \subset  \calW  \subset  \Sigma$ open sets.
  If $\calU_{\varphi(\alpha)} = (\upset_\Delta \varphi(\alpha) ) \setminus \aset{ \varphi(\alpha) }$ is connected 
  and $\inv{\varphi}(\calU_{\varphi(\alpha)}) \subset \calV$ for all $\alpha$ in $\calW \setminus \calV$,
  then any balanced map $m_\calV : \calV \to \ZZ_{\ge 1}$ extends to a balanced map $m_{\calW }: \calW \to \ZZ_{\ge 1}$ 
  by setting 
  \begin{align} 
    m_{\calW}(\alpha)  = \sum_{\substack{\gamma \in \inv \varphi(\beta) \\ \alpha \lessdot \gamma }} m_\calV(\gamma) 
  \end{align}
  for $\alpha \in \calW \setminus \calV$ and $\beta$ covering $\varphi(\alpha)$. 
  The value is independent of the choice of $\beta$.
\end{mytheorem}

While this result might seem unwieldy due to the conditions on $\calU_{\varphi(\alpha)}$, as a corollary we get the situation of interest for tropical moduli spaces.
Namely, if we got weights on top dimensional cells which satisfy the balancing condition, 
then they may be extended to an index map on the whole space $\Sigma$ whenever the target space $\Delta$ is \emph{strongly connected},
i.e.~$\Delta$~is connected and $\calU_\beta$ is connected for elements $\beta$ of rank at most $\dim \Delta - 2$.
 
\begin{myprop} 
    \label{myprop:ExtendingMap}
  Let $\Sigma$ and $\Delta$ be graded posets,
  $\dtmor : \Sigma \to \Delta$ a \combinatorial{} morphism,
  and $m_\calV$ balanced on codimension-1.
  If $\Delta$ is strongly connected, then $m_\calV$ extends to $m_\Sigma$.
\end{myprop}

    Finally, Lemma~\ref{lm:LiftingUpwardsPaths}  lifts all paths under certain conditions, which answers Question~\ref{que:Connectivity} with Proposition~\ref{prop:OneFibreConnection}.
    So we lift connectivity of tropical objects via \combinatorial{} maps.

\begin{myprop} 
    \label{myprop:LiftingConnectivity}
   Let $\Sigma$ and $\Delta$ be graded posets,
  $\dtmor : \Sigma \to \Delta$ a \combinatorial{} morphism,
  and $m_\Sigma$ a balanced map. 
  If $\Delta$ is connected in codimension-$k$, and $\varphi$ has a fibre that is connected in codimension-$k$ in $\Sigma$, 
  then $\Sigma$ is connected in codimension-$k$.
\end{myprop}

\subsection{Connections to past and future work}
    \label{sub:ConnectionsToPastAndFutureWork}
    Our three main influences are:
       the branched covers of cone complexes and of posets by Payne \cite{pay09};
       the approach to tropical moduli spaces by \cite {acp15, cmr16, ccuw20};
       and the combinatorial study of tropical morphisms undertaken in \cite{dv20}.    
    Our main motivation is the upcoming article \cite{var23}, which collects all the combinatorial constructions of \cite{dv20} which are balanced in codimension-1.
    By our work here, we get an indexed branched cover whose degree recovers the number of tropical morphisms that witness the gonality bound of a metric graph.
    
    For future research, we envision that this framework could give a different perspective to other works where tropical moduli spaces are constructed and covers are studied, e.g~\cite{cha12a, cmr16}.
    It is an open question how to extend our approach to the general case of spaces that are not necessarily polyhedral complexes. 
    So far, applications of this theory to tropical moduli spaces require passing to a barycentric subdivision, enormously complicating the combinatorics.
    A theory of poset stacks and their covers remains to be developed to better address the \emph{face category} associated with a tropical moduli space.

    \subsection{Acknowledgement}
The author thanks and acknowledges 
    Jan Draisma for illuminating discussions, mentorship throughout several years, and a careful reading of the first draft of this work; 
    Livio Ferretti for humouring questions on finite topological spaces; 
    and Martin Ulirsch for discussions on posets, category theory and notions of stacks. 
This research was supported by the Swiss National Science Foundation, grant number 200142.


\section{Polyhedral spaces}
        \label{sec:PolyhedralSpaces}


\mysubsection{Polyhedral spaces}
      \label{sub:PolyhedralSpaces}
We recall some facts of polyhedral geometry, mainly following \cite{zie12} and standard notation motivated by toric varieties, see \cite{cls11}.
Let $N$ be a free abelian group of finite rank, and write $N_\RR$ for the real vector space $N \otimes_\ZZ \RR$.
The elements of $N$ are the \emph{integral points} and those of  $N^* = \Hom(N, \ZZ)$ are the \emph{integral functionals}. 
To $u \in N^*$ and a constant $c \in \RR$ corresponds a rational closed upper half-space\index{half-space} $H^+(u,c) = \aset{x \in V \suchthat u(x) \ge c }$.
We consider rational polyhedra $\sigma$ in $N_\RR$, i.e~intersections of finitely many rational $H^+(u,c)$.
An integrally affine map $f : N_\RR \to N_\RR'$ is an affine map such that $f(N) \subset N'$; i.e.~$f(x) = y + L_\RR(x)$ with $y \in N'$ and $L_\RR$ 
the extension of some~$L \in \Hom(N, N')$.
These send rational polyhedra to rational polyhedra. 
The dimension $\dim \sigma$ of $\sigma$ is the dimension of its affine span $\affspan \sigma$.
A $k$-face of $\sigma$ is a face of dimension $k$, and a 0-face is a \emph{vertex}.
We denote by $N^\sigma$ the intersection $N \cap \affspan \sigma$.
The lineality space $L(\sigma)$ is the biggest subspace of $N_\RR$ with a translate contained in $\sigma$. 
We always assume trivial $L(\sigma)$.

\begin{de} 
We consider the following categories:
\begin{itemize} 
  \item $\polycat$ the category of pairs $(N, \sigma)$ of a finite-rank free abelian group $N$ and a polyhedron $\sigma$ in $N_\RR$ such that $L(\sigma) = \aset{0}$, with morphisms $f: (N, \sigma) \to (N', \sigma')$ given by integrally affine maps such that $f(\sigma) \subset \sigma'$.
  \item $\polyint$ the subcategory of $\polycat$ of pairs $(N, \sigma)$ such that $\sigma$ is rational, all the vertices of $\sigma$ are in $N$, and $N = N^\sigma$.
  \item $\polyintface$ with the same objects as $\polyint$ and morphisms $f : (\tau, N^\tau) \to (\sigma, N^\sigma)$ such that $N^\tau$ is mapped isomorphically to $N^\sigma \cap \affspan f(\tau)$ and $f(\tau)$ is a face of $\sigma$; these are called \emph{face morphisms}.
\end{itemize}
\end{de}

A face morphism $f: \sigma \to \sigma'$ is proper if $f(\sigma)$ is a proper face of~$\sigma'$. 
The non-proper face morphisms are the isomorphisms in $\polyintface$.
A polyhedral space is a collection of polyhedra glued along face morphisms. 

\begin{de} \label{de:PolyhedralSpace}
	Let $\Sigma$ be a finite category, and $\sigma : \Sigma \to \polyintface$ a functor $\alpha \mapsto (N^\alpha, \sigma_\alpha)$. 
  We slightly abuse notation and shorten $(N^\alpha, \sigma_\alpha)$ to just $\sigma_\alpha$ when it causes no confusion. 
  We say that $\sigma$ is a \emph{polyhedral space} if the following \emph{lifting conditions} are satisfied: 
	\begin{enumerate}[(a)]
    \item \label{item:WeGotAllFaces}
      For each $\alpha$ in $\Sigma$, and each proper face morphism $f_{\tau \sigma_\alpha}: (N, \tau) \to (N^\alpha, \sigma_\alpha)$ in $\polyintface$, there is a morphism $f$ in $\Sigma$ such that $\sigma(f) = f_{\tau \sigma_\alpha}$.
		
    \item \label{item:LiftingCondition}
      For any 
      $g : \gamma \to \alpha$ and 
      $h: \beta \to \alpha$ in $\Sigma$ 
      such that $\sigma(g)$ and $\sigma(h)$ are proper, 
      there is a bijection of diagrams in $\Sigma$ 
      to diagrams in $\polyintface$ of the form shown in Diagram~\ref{dia:fibred-category}.

	  \vspace{1em}	
		\begin{minipage}{0.9\linewidth}
			\centering
			
			\begin{tikzcd}
		& \beta \arrow[dr, "h"] & \\
		\alpha \arrow[rr, "g", swap] \arrow[ur, dashed, "f"] & & \gamma
		\end{tikzcd} \hspace{5em}
		\begin{tikzcd} 
		& \sigma_\beta \arrow[dr, "\sigma(h)"] & \\
		\sigma_\alpha \arrow[rr, "\sigma(g)", swap] \arrow[ur, dashed, "f_{\sigma_\alpha \sigma_\beta}"] & & \sigma_\gamma
		\end{tikzcd}
		
			\refstepcounter{thm}
			      \label{dia:fibred-category}
			Diagram~\ref{dia:fibred-category}
			\vspace{1em}
		\end{minipage}
	
		\noindent
		In other words, for every face morphism $f_{\sigma_\alpha \sigma_\beta} : \sigma_\alpha \to \sigma_\beta$ for which the diagram on the right commutes, 
		there is exactly one face morphism $f$ such that the diagram on the left commutes and $\sigma(f) = f_{\sigma_\alpha \sigma_\beta}$.
  \item \label{item:SkeletonCategory}
      The only isomorphisms in $\Sigma$ are self-maps.
	\end{enumerate}
\end{de}

\begin{re} 
    \label{re:PosetCategory}
  If $\sigma$ is a polyhedral space, 
  intuitively
  Condition~(a) implies that $\Sigma$ has non-redundant information on all the faces and containment of the polyhedra. 
  Condition~(b) is key to proving the stratification in Lemma~\ref{lm:stratification}.
  Condition~(c) is for the convenience of working with a skeleton category, as done in~\cite{acp15}, 
  so we get a \emph{face poset}, see Lemma~\ref{lm:FaceSetIsAPoset}.
\end{re}

\begin{re} 
  \label{re:PolyhedralSpaceConditions} 
  Definition~\ref{de:PolyhedralSpace} follows closely Definition 2.15 of \cite{ccuw20}, 
  but crucially we avoid having the lifting conditions for non-proper face morphisms, 
  because this would not be compatible with condition~(c). 
  Thus, we do not have a category fibered in groupoids over $\polyintface$, but this is not a problem since we do not pursue 2-categorical aspects.
\end{re}

\begin{de} 
   \label{de:MorphismPolyhedralSpace}
  A morphism $\Phi : [\sigma : \Sigma \to \polyintface] \to [\delta : \Delta \to \polyintface]$ is a pair $(\varphi, \, \aset{ \Phi_\alpha }_{\alpha \in \Sigma} )$
  of a functor $\varphi : \Sigma \to \Delta$ and a natural transformation $\aset{\Phi_\alpha : \sigma_\alpha \to \delta_{\varphi(\alpha)}}_{\alpha \in \Sigma}$ from $\sigma$ to $\delta \circ \varphi$ such that the image of $\Phi_\alpha$ is not contained in a proper face of $\delta_{\varphi(\alpha)}$. 

\begin{minipage}{0.95\textwidth}
	\centering
\[ 
\begin{tikzcd}
  \sigma_\alpha \arrow[r,"\sigma(f)"] 
\arrow[d,swap,"\Phi_\alpha"] &
\sigma_\beta  \arrow[d,"\Phi_\beta"] \\
 \delta_{\varphi(\alpha)} 
 \arrow[r,"\delta(\varphi(f))"] 
                           & \delta_{\varphi(\beta)}
\end{tikzcd} 
\]
	\refstepcounter{thm}
	      \label{dia:genrpc-morphism}
	Diagram~\ref{dia:genrpc-morphism}
	\vspace{1em}
\end{minipage} \\
That is, for every $f: \alpha \to \beta$ in $\Sigma$ we have that Diagram~\ref{dia:genrpc-morphism} in $\polyintface$ commutes. 
\end{de}

\begin{figure}
\centering

	\begin{overpic}[scale=0.8]{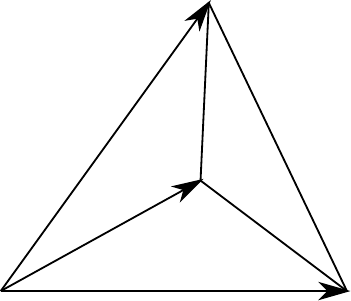}    
		\put (44,77) {\scalebox{0.8}{$\theta_1$}}
		\put (75,5) {\scalebox{0.8}{$\theta_2$}}
		\put (48,36) {\scalebox{0.8}{$\theta_3$}} 
	\end{overpic}
	\hspace{1em}
	\scalebox{0.8}{
		\begin{tikzpicture}
		\node (max) at (0,2) {$\aset{\theta_1, \theta_2, \theta_3 }$};
		\node (a) at (-2,1) {$\aset{\theta_1, \theta_2 }$};
		\node (b) at (0,1) {$\aset{\theta_1, \theta_3 }$};
		\node (c) at (2,1) {$\aset{\theta_2, \theta_3 }$};
		\node (d) at (-2,0) {$\aset{\theta_1 }$};
		\node (e) at (0,0) {$\aset{\theta_2 }$};
		\node (f) at (2,0) {$\aset{\theta_3 }$};
		\node (min) at (0,-1) {0};
		\draw (min) -- (d) -- (a) -- (max) -- (b) -- (f)
		(e) -- (min) -- (f) -- (c) -- (max)
		(d) -- (b);
		\draw[preaction={draw=white, -,line width=6pt}] (a) -- (e) -- (c);
		\end{tikzpicture}
	}
	\caption{\label{fig:simplicial} On the left, a simplicial cone $\sigma$ of dimension~3. On the right, the face diagram of $\sigma$ when embedded in $\PS$.}
\end{figure}
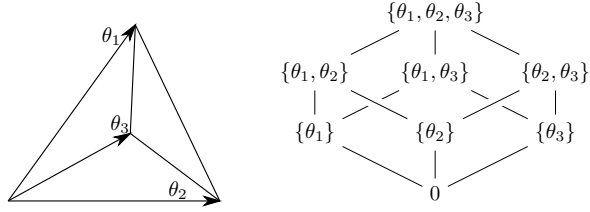

We denote by $\PS$ the category of polyhedral spaces.
When it causes no confusion we 
write $\Sigma$ for $\sigma : \Sigma \to \polyintface$.
Note that $\polyint$ embeds in $\PS$ by associating to $\sigma$ the polyhedral space of all faces of $\sigma$, e.g.~see Figure~\ref{fig:simplicial}.

\subsection{The face poset}
    \label{sub:TheFacePoset}
  Any category has a preorder on $\Obj \Sigma$ given by $\alpha \preceq \beta$ when $\Hom(\alpha, \beta) \ne \varnothing$.
  If $\Sigma$ indexes a polyhedral space we have

\begin{lm} 
      \label{lm:FaceSetIsAPoset}
  If $\sigma$ is a polyhedral space,
  then $(\Sigma, \preceq)$ is a poset called \emph{the face poset}.

\end{lm}

\begin{proof} 
    For antisymmetry, 
    if $\alpha \preceq \beta$ and $\beta \preceq \alpha$, 
    there are face morphisms $f: \alpha \to \beta$ and  $g : \beta \to \alpha$, 
    so $\dim \sigma_\alpha \le \dim \sigma_\beta$ and $\dim \sigma_\beta \le \dim \sigma_\alpha$.
    Hence, $\dim \sigma_\alpha = \dim \sigma_\beta$ and the face morphism $f$ is an isomorphism, 
    therefore $\alpha = \beta$ by Condition~(\ref{item:SkeletonCategory}) of Definition~\ref{de:PolyhedralSpace}.
\end{proof}

Given $\alpha$ in a poset $\Sigma$, the up-set and down-set generated by $\alpha$ are the posets 
\begin{align*}
\upset_{\Sigma} \alpha = \aset{ \gamma \in \Sigma \suchthat \alpha \le \gamma }, \quad \quad
\downset_{\Sigma} \alpha = \aset{ \gamma \in \Sigma \suchthat \alpha \ge \gamma }.
\end{align*} 
 
If the context allows, 
we write $\upset \alpha$ and $\downset \alpha$. 
The sets $\upset \alpha$ generate a topology, which we call the \emph{poset topology on} $\Sigma$, and whose open sets we also call up-sets.
The closed sets are generated by the $\downset \alpha$, and we also call them down-sets.

A morphism of posets is an order-preserving map $f: \Sigma \to \Delta$, i.e.~$x \le y$ implies~$f(x) \le f(y)$ for all $x,y \in \Sigma$.
This makes a category $\Poset$.          
The poset topology gives a fully faithful functor from $\Poset$ to $\Top$,
because a map of posets $f$ is continuous in the poset topology if and only if $f$ is order-preserving. 

\begin{re} 
    \label{re:FunctorIsOrderPreserving}
  A functor $\varphi : \Sigma \to \Delta$ is order-preserving for the preorder we have described in this subsection. 
  Thus, a morphism of polyhedral spaces induces a morphism of face posets, which is a continuous function in the poset topology.
\end{re}

\mysubsection{Topological realization}
        \label{sub:TopologicalRealization}
There is a faithful topological realization functor $\topo \cdot$ from  $\polycat$ to the category of topological spaces and continuous functions $\Top$. 
It maps $(N, \sigma)$ to $\sigma \subset N_\RR $ with the topology induced by $N_\RR$,
and morphisms $f$ are affine maps so automatically continuous.
As in Remark~2.2.1 of \cite{acp15}, 
all finite colimits exist in $\Top$,
so we get a faithful topological realization functor for $\PS$, 
sending $\Sigma$ to the colimit 
\[ \topo \Sigma = \colim_{\alpha \in \Sigma} \topo{\sigma_\alpha} \in \Top.\] 
So we denote by $p_\alpha : \topo{\sigma_\alpha} \to \topo{\Sigma}$ the universal maps
satisfying $p_\alpha = p_\beta \circ \sigma(f)$ for any $f : \alpha \to \beta$ in $\Sigma$,
and $\topo \Sigma$ carries the final topology with respect to all the $p_\alpha$. 
For $\Phi : \Sigma \to \Delta$,
the topological realization of $\Phi$ is a continuous function $\topo \Phi : \topo \Sigma \to \topo \Delta$ obtained by gluing all the $\topo {\Phi_\alpha}$, 
such that for all $\alpha \in \Sigma$ the Diagram~\ref{dia:HowTopologicalRealizationGlues} commutes. 

  \begin{minipage}{0.95\textwidth}
	\centering
\[ 
\begin{tikzcd}
    \topo {\sigma_\alpha} \arrow[r,"p_\alpha"] \arrow[d,swap,"\topo {\Phi_\alpha}"] & \topo \Sigma \arrow[d,"\topo \Phi"] \\
	\topo {\delta_{\varphi(\alpha)}} \arrow[r,"p_{\varphi(\alpha)}"] & \topo \Delta
\end{tikzcd} 
\]
	\refstepcounter{thm}
	      \label{dia:HowTopologicalRealizationGlues}
	Diagram~\ref{dia:HowTopologicalRealizationGlues}
	\vspace{1em}
\end{minipage}


\mysubsection{Relative interior}
    \label{RelativeInteriorOfSigma}
Let $\sigma \subset N_\RR$ be a polyhedron with $\affspan \sigma$ not necessarily equal to~$N_\RR$.
The \emph{relative interior} $\sigma^\circ$ of $\sigma$ is the topological interior of $\topo \sigma$ as a subspace of $\affspan \sigma$.
This equals the complement in $\topo \sigma$ of the union of all the proper faces of~$\sigma$;
see e.g.~\cite[Section 2.3]{zie12} for a reference.
The relative interior commutes with affine maps, namely~$f(\sigma^\circ) = f(\sigma)^\circ$. 
Moreover, for maps in $\polycat$ we have

\begin{lm} 
    \label{lm:AffineMapPreservesInterior} 
    Let $f : (N^\sigma, \sigma) \to (N^\delta, \delta)$ be an affine map in $\polycat$. 
  If $f(\sigma)$ is not contained in a proper face of $\delta$, then $f(\sigma^\circ) \subseteq \delta^\circ$.
  Equivalently, if $f(\sigma^\circ) \cap \delta^\circ \ne \varnothing$, then $f(\sigma^\circ) \subseteq \delta^\circ$.
\end{lm}

\begin{proof}
  This fact follows from $f$ being an affine map, so $f(\sigma)$ is a polyhedron in $\affspan N_\delta$, and from $f(\sigma)^\circ = f(\sigma^\circ)$.
  See e.g.~\cite[Lemma 9.13]{var22} for details.
\end{proof}

The relevance of Lemma~\ref{lm:AffineMapPreservesInterior} is that the maps $\Phi_\alpha$ from a morphism $(\varphi, \aset{\Phi_\alpha}_{\alpha \in \Sigma})$ of polyhedral spaces satisfy it. 
With a face morphism we get a stronger result.
\begin{lm}  
  \label{lm:IntersectionIffMorphism}
Let $\Sigma \in \PS$ and
$\alpha, \beta$ in $\Obj \Sigma$.
There is a morphism $\alpha \to \beta$ if and only if $p_\alpha(\sigma_\alpha^\circ) \cap \im p_\beta \ne \varnothing$.
Moreover, if $p_\alpha(\sigma_\alpha^\circ) \cap \im p_\beta \ne \varnothing$, then $p_\alpha(\sigma_\alpha^\circ) \cap \im p_\beta = p_\alpha(\sigma_\alpha^\circ)$.
\end{lm}
\begin{proof}
  For convenience, we write $\sigma f$ for $\sigma(f)$. 
  Let $f_{\alpha\beta}: \alpha \to \beta$ be a morphism. 
  By the universal property of $\topo \Sigma =  \colim_{\alpha \in \Sigma} \topo{\sigma_\alpha}$ we  get $p_\alpha = p_\beta \circ \sigma f_{\alpha\beta}$.
  Thus, $\im p_\alpha \subset \im p_\beta$.

  On the other direction, suppose there is $z \in p_\alpha(\sigma_\alpha^\circ) \cap \im p_\beta $. 
  Let $x \in { \sigma^\circ_\alpha}$ and $y \in {\sigma_{\beta}}$ such that $p_\alpha(x) = p_{\beta}(y) = z$,
  and let $\tau$ be the unique face of $\sigma_\beta$ such that $y \in \tau^\circ$. 
  By Condition~(a) of Definition~\ref{de:PolyhedralSpace} there is $\alpha'$ and $\iota : \alpha' \to \beta$ in $\Sigma$ 
  such that $\sigma \iota (\sigma_{\alpha'}) = \tau \subset \sigma_\beta$.
  Since  $\sigma \iota$ is bijective, let $x'$ in $\sigma_{\alpha'}$ be the preimage of $y$ under $\sigma \iota$.
  Note that $x'$ is in $\sigma_{\alpha'}^\circ$.
  
  Since $p_{\alpha'}(x') = z = p_\alpha(x)$, 
  by the universal property of $\topo \Sigma$ there is a sequence of face morphisms connecting $x$ and $x'$; 
  e.g.~$\alpha \xrightarrow{\, f_1 } \alpha_1 \xleftarrow{\, f_2 } \alpha_2 \xleftarrow{\, f_3 } \dots \xrightarrow{\, f_n } \alpha'$ with $x' =  \sigma f_n(\dots \sigma \inv  f_3(\sigma \inv f_2( \sigma f_1(x) ) ) )$. 
  We claim one can replace the first two morphisms by a single one,
  so iterating, one gets a length-1 sequence $f : \alpha \to \alpha'$ with $\sigma f(x) = x'$. 
  If $\alpha \ne \alpha'$, then Condition~(c) of Definition~\ref{de:PolyhedralSpace} implies that $\dim \alpha \ne \dim \alpha'$, 
  but then $\sigma f (x)$ is contained in a proper face of $\sigma_{\alpha'}$, contradicting that $x' = \sigma f (x)$ is in $ \sigma_{\alpha'}^\circ$. 
  So $\alpha = \alpha'$, and $\iota$ is the sought morphism $\alpha \to \beta$.
  
  To prove the claim,
  there are 4 cases for the first 2 morphisms:
	$\alpha \xrightarrow{\, f_1 } \alpha_1 \xrightarrow{\, f_2 } \alpha_2$, 
	$\alpha \xleftarrow{\, f_1 } \alpha_1 \xleftarrow{\, f_2 } \alpha_2$, 
	$\alpha \xleftarrow{\, f_1 } \alpha_1 \xrightarrow{\, f_2 } \alpha_2$, and 
	$\alpha \xrightarrow{\, f_1 } \alpha_1 \xleftarrow{\, f_2 } \alpha_2$. 
	In the first two we compose them. 
  In the third, $\sigma f_1(\sigma_{\alpha_1})$ is a face of $\sigma_\alpha$ that contains the interior point $x$. 
  Hence, $f_1$ is an isomorphism, and we can reverse the arrow to fall into the first case. 

  \vspace{1em}
	\begin{minipage}{0.95\linewidth}
		\centering
		
		\begin{tikzcd}
		& \alpha_2 \arrow[dr, "f_2"] & \\
		\alpha \arrow[rr, "f_1", swap] \arrow[ur, dashed, "\tilde f"] & & \alpha_1
		\end{tikzcd} \hspace{5em}
		\begin{tikzcd} 
		& \sigma_{\alpha_2} \arrow[dr, "\sigma f_2"] & \\
		\sigma_{\alpha} \arrow[rr, "\sigma f_1", swap] \arrow[ur, dashed, ""] & & \sigma_{\alpha_1}
		\end{tikzcd}
		
		\refstepcounter{thm}
		      \label{dia:fourth-case}
		Diagram~\ref{dia:fourth-case}
		\vspace{1em}
	\end{minipage}

	The fourth gives Diagram~\ref{dia:fourth-case}. 
  There are points $x_2 \in \sigma_{\alpha_2}$ and $x_1 \in \sigma_{\alpha_1}$, 
  and face morphisms $f_1 : N^{\alpha} \to N^{\alpha_1}$ and $f_2 : N^{\alpha_2} \to N^{\alpha_1}$,
  such that $x_1 = \sigma f_1(x) = \sigma f_2(x_2)$.
  Both $\im \sigma f_1$ and $\im \sigma f_2$ are faces of $\sigma_{\alpha_1}$ that contain $x_1$.
  Since $x$ is in $\sigma_\alpha^\circ$, the point $x_1 = \sigma f_1(x)$ is an interior point of $\im \sigma f_1$, thus the face $\im \sigma f_1$ is contained in $\im \sigma f_2$. 
  As both maps are face morphisms, they are injective and affine, so $\inv{(\sigma f_2)} \circ \sigma f_1$ restricts to a face morphism $\sigma_\alpha \to \sigma_{\alpha_2}$. 
  %
  By Condition~(b) of Definition~\ref{de:PolyhedralSpace}, there exists a unique lift $\tilde f: \alpha \to \alpha_2$, 
  and this replaces $\alpha \xrightarrow{\, f_1 } \alpha_1 \xleftarrow{\, f_2 } \alpha_2$.
	  \qedhere 
 \end{proof} 


\mysubsection{Decomposing the topological realization}
      \label{sub:DecomposingTheTopologicalRealization}
From Lemma~\ref{lm:IntersectionIffMorphism} we get a stratification of $\topo \Sigma$ as a disjoint union of interiors, analogous to embedded polyhedral complexes.



\begin{lm} 
    \label{lm:stratification}
	Let $\sigma : \Sigma \to \polyintface$ be in $\PS$. 
  We have that 
	\[ \topo \Sigma  = \bigsqcup_{\alpha \in \Sigma} p_\alpha(\sigma_{\alpha}^\circ). \]
\end{lm}

\begin{proof}
  By the universal property of $\topo \Sigma$, for every $x$ in $\topo \Sigma$ there is $\gamma$ in $\Sigma$ with $x \in \im p_\gamma$.
  Choose $\hat x$ in the fibre $\inv p_\gamma(x)$, and $\tau$ a face of $\sigma_\gamma$ with $\hat x$ in $\tau^\circ$. 
  By Property~(a) of Definition~\ref{de:PolyhedralSpace} there is $\alpha$ and $f: \alpha \to \gamma$ in $\Sigma$ such that $\sigma f(\sigma_\alpha) = \tau$.
  Thus, $\hat x$ is in $\tau^\circ = \sigma f(\sigma_\alpha)^\circ = \sigma f(\sigma_\alpha^\circ)$,
  and $x = p_\gamma(\hat x)$ is in $(p_\gamma \circ \sigma f )(\sigma_\alpha^\circ) = p_\alpha(\sigma_\alpha^\circ)$,
so indeed~$\topo \Sigma = \bigcup p_\alpha(\alpha^\circ_\alpha)$. 
  By Lemma~\ref{lm:IntersectionIffMorphism} and Property~(c) of Definition~\ref{de:PolyhedralSpace} this union is disjoint.  \qedhere
\end{proof}

Thus, we can describe intersections $\im p_\alpha \cap \im p_\beta$.
In the following, $\downset \alpha$ is as in Subsection~\ref{sub:TheFacePoset}, and equals the set of domains of the maps in $\Hom(-, \alpha)$.

\begin{lm}
    \label{lm:impalphaIntersectsimpbeta}
  Let $\sigma : \Sigma \to \polyintface$ be in $\PS$, and $\alpha, \beta$ in $\Sigma$.
  We have that 
  \[ \im p_\alpha \cap \im p_\beta = \bigcup_{\eta \in \downset \alpha \cap \downset \beta} \im p_\eta. \]
\end{lm}

\begin{proof} 
  If $x$ is in $\im p_\alpha \cap \im p_\beta$, then Lemma~\ref{lm:stratification} gives an $\eta \in \Sigma$ with $x \in p_\eta(\sigma^\circ_\eta)$, so Lemma~\ref{lm:IntersectionIffMorphism} gives morphisms $\eta \to \alpha$ and $\eta \to \beta$. 
  Conversely, if there are morphisms $f:\eta \to \alpha$ and $g:\eta \to \beta$, we have $p_\eta = p_\alpha \circ \sigma f$  and $p_\eta = p_\beta \circ \sigma g$, so $\im p_\eta \subset \im p_\alpha \cap \im p_\beta$.
\end{proof}

Note that $\Hom_\Sigma(\alpha, \alpha)$ equals $\Aut \alpha$, thus
$p_\alpha$ factors through the projection to the quotient space $\topo {\sigma_\alpha} / \Aut \alpha$ of the action of $\Aut \alpha$ on $\topo {\sigma_\alpha}$.  
The map $\overline p_\alpha :  \topo {\sigma_\alpha} / \Aut \alpha \to \topo \Sigma $ is an homeomorphism onto its image; see Diagram~\ref{dia:orbifold}. Thus, by Lemma~\ref{lm:stratification}, $\topo \Sigma$ decomposes as a disjoint union of interiors of polyhedra modulo automorphisms. 

\begin{minipage}{0.95\textwidth}
	\center
	\[ 
	\begin{tikzcd}[row sep=tiny]
	\lvert \sigma_\alpha \rvert \arrow[dd,swap] \arrow[dr,"p_\alpha"]&  \\
	& \lvert\Sigma\rvert\\
	\lvert \sigma_\alpha \rvert / \Aut \alpha \arrow[ur, swap, "\overline p_\alpha"] & 
	\end{tikzcd} 
	\]
	
	\refstepcounter{thm}
	      \label{dia:orbifold}
	Diagram~\ref{dia:orbifold}
	\vspace{1em}
\end{minipage}

\mysubsection{Polyhedral complexes and barycentric subdivision} 
      \label{sub:Polyhedralcomplexes}
      If $(\Obj \Sigma, \preceq)$ encodes all the information of $\Sigma$, 
      we get trivial $\Aut \alpha$ for all $\alpha \in \Sigma$.
      Thus, the $p_\alpha$ in Diagram~\ref{dia:orbifold} are injective, 
      and $\topo \Sigma$ is locally polyhedral in a straightforward manner.  

\begin{de} 
      \label{de:PolyhedralComplex}
    A polyhedral space $\sigma : \Sigma \to \polyintface$ is a 
    \emph{rational polyhedral complex}
    if $\card{\Hom_\Sigma(\alpha, \beta)} \le 1$.
    We write $f_{\alpha\beta}$ for the unique element in $\Hom_\Sigma(\alpha, \beta)$, if there is one,
    and $\RPC$ for the full subcategory of $\PS$ of polyhedral complexes.
\end{de}

\begin{ex}[embedded polyhedral complex]
      \label{ex:EmbeddedPolyhedralComplex}
      Let $(V,N)$ be a vector space with an integral structure, 
      and $\Sigma$ a finite family of integral polyhedra in $(V,N)$ such that:
      \begin{enumerate}[(a)]
           \item if $\sigma$ is in $\Sigma$ and $\tau$ is a face of $\sigma$, then $\tau$ is in $\Sigma$;
           \item if $\sigma_1, \sigma_2$ are in $\Sigma$ and $\sigma_1 \cap \sigma_2 \ne \varnothing$,
             then $\sigma_1 \cap \sigma_2$ is a face of $\sigma_1$ and a face of $\sigma_2$.
      \end{enumerate}
      The family $\Sigma$ is a category with set inclusions as morphisms, 
      so with $\sigma$ the identity functor,
      we get a polyhedral complex per Definition~\ref{de:PolyhedralSpace}. 
      We have that 
      \begin{align} 
           \label{eq:TopologicalRealizationEmbeddedPolyComplex} 
           \topo \Sigma \sim \bigcup_{\sigma \in \Sigma} \sigma.
      \end{align}
      So topological realization generalizes the support of an embedded polyhedral complex.
\end{ex}

While non-trivial groups of automorphisms frequently arise in the construction of moduli spaces, 
focusing on polyhedral complex is justified in the context of tropical moduli spaces by refining to the barycentric subdivision, 
at the cost of complicating the combinatorics of the face poset.

\begin{de}
      \label{de:Refinement}
  A \emph{refinement} of a polyhedral space $\Sigma \in \PS$ is a morphism $\Phi : \tildeSigma \to \Sigma \in \PS$ 
  such that the topological realization $\topo \Phi : \topo {\Sigma'} \to \topo \Sigma$ 
  is a homeomorphism and induces a bijection on integral points, i.e.~$\topo \Phi(N') = N$. 
\end{de}

\begin{ex} 
    \label{ex:HalfSpace}
    Given $\sigma \subset N_\RR$, and hyperplane $H$ such that $H \cap \sigma^\circ \ne \varnothing$ induces a non trivial refinement of $\sigma$, namely $\sigma \cap H^-$ and $\sigma \cap H^+$. Moreover, if $\sigma$ has unbounded faces, $H$ can be chosen such that all the bounded faces are contained in $H^-$, so $\sigma \cap H^-$ is bounded and $\sigma \cap H^+$ has a unique bounded face.
\end{ex}

We sketch the barycentric subdivision for polyhedral spaces whose cells $\sigma \subset (V, N)$ are cones.
Denote by $\coneintface$ the subcategory of $\polyintface$ whose objects are cones.
Recall that the polyhedra in $\polyintface$ have trivial lineality space, hence the cones of $\coneintface$ are pointed, i.e.~they have a 0-face and are spanned by their 1-faces.
If $\sigma$ is rational and $0 \lessdot \theta \preceq \sigma$ is a 1-face, 
$\theta \cap N$ is a monoid generated by the \emph{primitive vector} $v_\theta$.
Face morphisms $f$ respect primitive vectors, i.e.~$f(v_\theta) = v_{f(\theta)}$. 
Thus if we define the \emph{barycentre} of $\sigma$ as
\[\barycentre \sigma = \sum_{0 \lessdot \theta \preceq \sigma} v_\theta,\]
we get that $f(\barycentre \sigma) = \barycentre{f(\sigma)}$.
With this we construct a refinement:
 
\begin{cons}[Barycentric Subdivision] 
    \label{cons:BarycentricSubdivision}
   Let $\sigma : \Sigma \to \coneintface$ be a polyhedral space of cones.
   Denote by $\BCS{\Sigma}$ the set of all chains $\calC : \alpha_1 \xrightarrow{f_1} \alpha_2  \xrightarrow{f_2} \dots  \xrightarrow{f_{k}} \alpha_k$ of proper morphisms in $\Sigma$.
   To $\calC \in \BCS{\Sigma}$ we associate the rational cone
   \[ \BCS{\sigma}(\calC) = (\Span_{\RRgo}(\barycentre{\sigma(\alpha_k)}, \barycentre{\sigma f_{k-1}(\sigma(\alpha_{k-1}))}, \dots, \sigma (f_{k-1} \circ \dots \circ f_{1})(\sigma(\alpha_{1})) ), N^{\alpha_k}). \]
   To get a functor $\BCS{\sigma} : \BCS \Sigma \to \coneintface$ that is a polyhedral space, 
   we need to add  to $\BCS \Sigma$ an object to act as the vertex,
   morphisms to satisfy the lifting conditions, 
   and afterwards to pass down to a skeleton category.
   For every pair of chains $\calC : \dots \xrightarrow{f_{k}} \alpha_k$ and $\calD : \dots \xrightarrow{g_{l}} \gamma_l$ and every morphism $h : \alpha_k \to \gamma_l$,
   add a morphism $h^{\calC, \calD} : \calC \to \calD$. 
   The associated polyhedral map $\BCS{\sigma}(h^{\calC,\calD})$ is given by restricting $\sigma(h)$ to $\BCS{\sigma}(\calC)$,
   which maps to $\BCS{\sigma}(\calD)$ because $f(\barycentre{\sigma}) = \barycentre{f(\sigma)}$.
   Finally $h_2^{\calD, \calE} \circ h_1^{\calC, \calD} = h_3^{\calC, \calE}$ if and only if there are $h_1, h_2, h_3 \in \Mor \Sigma$ such that $h_2 \circ h_1 = h_3$.
\end{cons}

Note that for $\calC : \dots \xrightarrow{f_{k}} \alpha_k$ we have $\dim  \BCS{\sigma}(\calC) = k+1$ and that the cone is simplicial, i.e.~isomorphic to $\RR^{k+1}_{\ge 0}$.
Moreover, as pointed out in \cite{acp15}, this construction resolves away automorphisms and multiple glueings of one face to one cone. 

\begin{prop} 
    \label{prop:BCSIsPolyhedralComplex}
    Let $\sigma : \Sigma \to \coneintface$ be a \pscone{}.
    The barycentric subdivision $\BCS \sigma$ is a polyhedral complex of simplicial cones.
\end{prop}

\begin{re} 
    \label{re:GeneralizeBarycentricSubdivision}
    A compact polyhedron $\sigma \subset (V, N)$ is called a \emph{polytope}. 
    To define its barycentre we consider the \emph{homogenization} $\tilde \sigma$ of $\sigma$, i.e.~embed $(V,N)$ into the plane of $(V \times \RR, N \times \ZZ )$ whose last coordinate equals~1 and consider $\tilde \sigma =  \Span_{\RRgo} \sigma$.
    Define $\beta(\sigma)$ to be equal to $\sigma \cap (\Span_{\RRgo} \beta(\tilde \sigma))$.
    Similarly, the barycentric subdivision $\BCS{\tilde \sigma}$ of $\tilde \sigma$ induces a barycentric subdivision $\BCS \sigma$ on $\sigma$.
    Writing $d$ for $\dim \sigma$, the cells of $\BCS \sigma$ are isomorphic to the $d$-simplex $\Delta^d$, i.e.~the convex hull $\operatorname{convex} \aset{e_1, \dots, e_{d+1}}$ in $\RR^{d+1}$ of the $d+1$ standard vectors.

    Any polyhedron $\sigma$ can be decomposed as a Minkowski sum of a polytope $P(\sigma)$ and a cone $C(\sigma)$.
    The cone $C(\sigma)$ is called the \emph{recession cone} of $\sigma$.
    Thus, $\BCS{P(\sigma)}$ and $\BCS{C(\sigma)}$ induce a refinement of $\sigma$ such that each cell is a Minkowski sum of a simplex and a simplicial cone.
    This construction resolves away automorphisms of any polyhedron, thus generalizing Proposition~\ref{prop:BCSIsPolyhedralComplex}.
\end{re}

\mysubsection{Relating the topology of \texorpdfstring{$(\Sigma, \preceq)$}{(Σ, ⪯)} and \texorpdfstring{$\topo \Sigma$}{|Σ|}}
      \label{sub:RelatingThePosetAndTheFinalTopology}
      Given $\Sigma \in \PS$ and $x \in \topo \Sigma$, 
      consider the set $\aset{ \gamma \in \Sigma \suchthat x \in \im p_\gamma }$. 
By Lemmas~\ref{lm:IntersectionIffMorphism} and~\ref{lm:stratification} this is an up-set generated by the element $\alpha$ in $\Sigma$ such that $x \in p_\alpha(\sigma_\alpha^\circ)$.
So the following map is well defined: 

\begin{de}
    \label{de:PolyFunction}
  Let $\Sigma$ be in $\PS$.	
  We define the map $\poly_\Sigma : \topo \Sigma \to \Sigma$ given by 
  \[ x \mapsto  \min \, \aset{ \gamma \in (\Sigma, \preceq) \suchthat x \in \im p_\gamma }. \]
\end{de}
Now we assume that $\Sigma$ is a polyhedral complex, so the universal maps $p_\alpha : \sigma_\alpha \to \topo \Sigma$ have inverses. 
 
\begin{lm}
  \label{lm:PolyIsSurjectiveOpenContinuous}
	Let $\sigma: \Sigma \to \polyintface $ be in $\RPC$.
	The map $\poly_\Sigma : \topo \Sigma \to \Sigma$ is surjective, open, and continuous.
\end{lm}

\begin{proof}
	Since $\poly_\Sigma(x)$ equals $\alpha$ for any $x \in p_\alpha(\sigma_\alpha^\circ)$, we have surjectivity.
	
	Recall that $V \subseteq \topo \Sigma$ is open in the final topology on $\topo \Sigma$ if and only if $\inv p_\alpha(V)$ is open in $\topo {\sigma_\alpha}$ for all $\alpha \in \Sigma$. 
	To see that $\poly_\Sigma$ is open, let $V \subset \topo \Sigma$ be an open set, $\alpha \in \poly_\Sigma(V)$, and $\beta$ in $\Sigma$ such that $\alpha \preceq \beta$. 
  There is $x$ in $V$ such that $x \in p_\alpha(\sigma_\alpha^\circ )$. 
  Since $\alpha \preceq \beta$, there is $f_{\alpha\beta} : \alpha \to \beta \in \Sigma$. 
  Thus, $p_\alpha = p_\beta \circ \sigma(f_{\alpha\beta})$
  and we have $x \in \im p_\beta$, 
  so $\inv p_\beta(V)$ is an open neighbourhood of $\inv p_\beta(x)$ in $\topo {\sigma_\beta}$. 
  There is a point $y$ in $\inv p_\beta(V) \cap \sigma_\beta^\circ$, because any neighbourhood of a point in $\topo {\sigma_\beta}$ intersects $\sigma_\beta^\circ$.
  This gives that $p_\beta(y) \in V$ and $\poly_\Sigma (p_\beta(y)) = \beta$.
  Hence $\poly_\Sigma(V)$ is an up-set, namely open.
	
	Finally, $\aset{\upset \alpha \suchthat \alpha \in \Sigma}$ generates the poset topology, so  $\poly_\Sigma$ is continuous if 
  \begin{align} 
    \label{eq:continuityOpen}
   \inv {\poly_\Sigma}(\upset \alpha) = 
  \bigcup_{\alpha \preceq \gamma} p_\gamma(\sigma_\gamma^\circ)
  \end{align}
  is open, 
  i.e.~$\inv p_\beta(\inv {\poly_\Sigma}(\upset \alpha)) \subset N^\beta_\RR$ is open for all $\alpha, \beta$ in $\Sigma$. 
  By Lemma~\ref{lm:IntersectionIffMorphism} we have that $p_\beta(\sigma_\beta) \cap p_\gamma(\sigma_\gamma^\circ) \ne \varnothing$ if and only if $\gamma \preceq \beta$. 
  In that case, 
  $p_\beta(\sigma_\beta) \cap p_\gamma(\sigma_\gamma^\circ) = p_\gamma(\sigma_\gamma^\circ)$, so: 
  \begin{align*} 
    \inv {\poly_\Sigma}(\upset \alpha) \cap  p_\beta(\sigma_\beta) &= 
  \bigcup_{\alpha \preceq \gamma} p_\gamma(\sigma_\gamma^\circ) \cap p_\beta(\sigma_\beta) 
  = \bigcup_{\alpha \preceq \gamma \preceq \beta} p_\gamma(\sigma_\gamma^\circ) \\
  &= \bigcup_{\alpha \preceq \gamma \preceq \beta} p_\beta \circ \sigma(f_{\gamma\beta})(\sigma_\gamma^\circ) 
  = p_\beta ( \bigcup_{\alpha \preceq \gamma \preceq \beta} \sigma(f_{\gamma\beta})(\sigma_\gamma^\circ)).
  \end{align*}
Applying $\inv p_\beta$ on both sides we get
	\[ 
    \inv p_\beta(\inv {\poly_\Sigma}(\upset \alpha)) =  
  \bigcup_{\alpha \preceq \gamma \preceq \beta} \sigma(f_{\gamma \beta})(\sigma_{\gamma}^\circ), 
\]
and it is straightforward to show this is a complement of faces of $\sigma_\beta$, thus an open set, i.e.
\[
\topo {\sigma_\beta}=  
  \bigcup_{\alpha \preceq \gamma \preceq \beta} \sigma(f_{\gamma \beta})(\sigma_{\gamma}^\circ)
  \sqcup
  \bigcup_{\substack{\eta \preceq \beta \\ \alpha \not \preceq \eta}} \sigma(f_{\eta \beta})(\sigma_\eta).  \qedhere
\]
\end{proof}

\begin{lm}
  \label{lm:superDiagram}
	Let $\Phi : \Sigma \to \Delta$ be a morphism in $\RPC$, and $\topo \Phi : \topo \Sigma \to \topo \Delta$ the induced map on the topological realizations. 
  We have that Diagram~\ref{dia:super-diagram} commutes.
	
	\begin{minipage}{0.95\textwidth}
	\centering
	\vspace{0.7em}
	\[ 
	\begin{tikzcd}
	\topo \Sigma \arrow[r,"\poly_\Sigma"] \arrow[d,swap,"\topo \Phi"] &
	\Sigma \arrow[d,"\varphi"] \\
	\topo \Delta \arrow[r,"\poly_\Delta"] & \Delta
	\end{tikzcd} 
	\]
	
	\refstepcounter{thm}
       	\label{dia:super-diagram}
	Diagram~\ref{dia:super-diagram}
	\vspace{1em}
	\end{minipage}
\end{lm}

\begin{proof}
  Let $\Phi  = (\varphi, \aset{\Phi_\alpha }_{\alpha \in \Sigma} )$,
	a point $x$ in $\topo \Sigma$, and $\alpha = \poly_\Sigma(x)$.
	So $\varphi \circ \poly_\Sigma(x) = \varphi(\alpha)$.
	On the other hand, by Lemma~\ref{lm:AffineMapPreservesInterior}, and the discussion succeeding it, we have that $\Phi_\alpha(\sigma_\alpha^\circ) \subseteq \delta^\circ_{\varphi(\alpha)}$.
	Thus,  $\Phi_\alpha( \inv p_\alpha (x) )$ is in~$\delta_{\varphi(\alpha)}^\circ$, therefore $\poly_\Delta \circ p_{\varphi(\alpha)} \circ \Phi_\alpha \circ \inv p_\alpha(x) = \varphi(\alpha)$.
	By Diagram~\ref{dia:HowTopologicalRealizationGlues}  we have $p_{\varphi(\alpha)} \circ \Phi_\alpha \circ \inv p_\alpha(x) =\topo \Phi(x)$, so we are done.
\end{proof}

From Diagram~\ref{dia:super-diagram} we get a containment
\begin{align} 
  \label{eq:containment-fibres}
  \poly_\Sigma (\inv {\topo \Phi} (y)) \subset \inv \varphi( \poly_\Delta(y) ).
\end{align}
So one can wonder if $\poly_\Sigma$ maps the fibre $ \inv {\topo \Phi} (y)$  bijectively onto the fibre $ \inv \varphi( \poly_\Delta(y) )$, 
which is desirable for counting purposes.
Clearly not: consider embedded polyhedral complexes $\topo \Sigma, \topo \Delta$, and $\topo \Phi$ a proper containment relation $\topo \Sigma \subset \topo \Delta$.
Any point $y \in \topo \Delta \setminus \topo \Sigma$ disproves the bijection, but under a technical condition on $\dtmor$ 
we do get it in in Lemma~\ref{lm:CombinatorialImpliesFibreBijection}.

 
           \section{Indexed branched covers} 
                \label{sec:IndexedBranchedCovers}


\mysubsection{Indexed branched covers} 
        \label{sub:IndexedBranchedCovers}
Inspired by maps of Riemann surfaces, see Example~\ref{ex:LocalDegreeRealizability}, 
for us a pair $(F, B)$ of a continuous map $F : X \to Y$ and a closed subset $B \subset Y$,
is a \emph{branched cover}
if $Y \setminus B$ is open and dense, and $F$ restricted to $\inv F (Y \setminus B)$ is a cover.  
The set $B$ is called the \emph{branch locus}.

\begin{figure}
\centering
	\begin{minipage}{0.24\textwidth}
		\scalebox{0.9}{
		\begin{overpic}{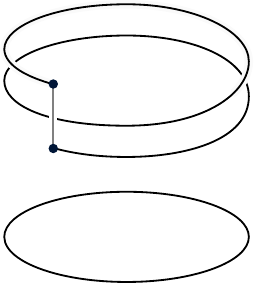} 
			\put (96,38) {\scalebox{1.3}{$\downarrow$}}
			
			\put (95,60) {\scalebox{1}{$X$}} 
			\put (95,13) {\scalebox{1}{$Y$}} 
			
		\end{overpic}
	}
	\end{minipage} 
\begin{minipage}{0.37\textwidth}
	\scalebox{0.8}{
		 \begin{tikzpicture}
		\node (A1) at (0,0) {$A_1$};
		\node (B1) at (2,0) {$B_1$};
		\node (C1) at (3.6,0) {$C_1$};
		\node (C2) at (4.4,0) {$C_2$};
		\node (t1) at (2.5,1) {$t_1$};
		\node (t2) at (3.5,1) {$t_2$};
		\node (s1) at (0.5,1) {$s_1$};
		\node (s2) at (1.5,1) {$s_2$};
		\draw (A1) -- (s1) -- (B1) -- (t1) -- (C1);  
		\draw (A1) -- (s2) -- (B1) -- (t2) -- (C2);  
		 \draw[preaction={draw=white, -,line width=6pt}] (B1) -- (s1);
		 \draw[preaction={draw=white, -,line width=6pt}] (B1) -- (t2);
		 \node (A) at (0,-2) {$A$};
		 \node (B) at (2,-2) {$B$};
		 \node (C) at (4,-2) {$C$};
		 \node (t) at (3,-1) {$t$};
		 \node (s) at (1,-1) {$s$};
		 \draw (A) -- (s) -- (B) -- (t) -- (C);  
		 \node (t2) at (5,0.1) {\scalebox{1.25}{$\Gamma$}};
		 \node (t2) at (5,-0.8) {\scalebox{1.4}{$\downarrow$}};
		 \node (t) at (5,-1.7) {\scalebox{1.25}{$\Delta$}};
		 \end{tikzpicture}
	}
\end{minipage}
\begin{minipage}{0.25\textwidth}
	\begin{overpic}{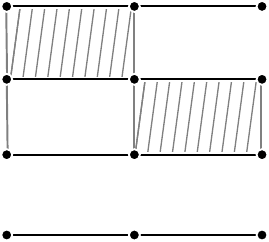} 
		\put (110,30) {\scalebox{1.3}{$\downarrow$}}
		
		\put (107,55) {\scalebox{1}{$\topo \Gamma$}} 
		\put (105,-2) {\scalebox{1}{$\topo \Delta$}} 
		
	\end{overpic}
\end{minipage}
	\caption{\label{fig:covers} 
  On the left, a degree-2 cover of a circle, where the vertical grey line indicates identification of two points. 
  In the centre,  a branched covering of posets.  
  On the right, a branched covering of a segment subdivided in two pieces.  
  Both branched coverings admit an index map that makes them a degree-3 indexed branched covering.}
\end{figure}
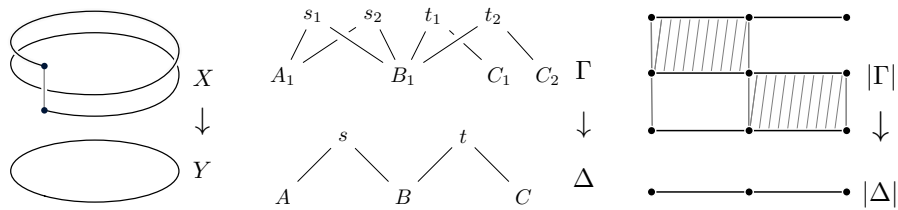

\begin{ex} 
      \label{ex:BranchedCoverPosets}
  Let $\varphi : \Sigma \to \Delta$ be a morphism of posets. 
  The closure of a set $\calU \subset \Delta$ is $\downset_\Delta \calU$, 
  so a subset of $\Delta$ is dense if and only if it contains~$\max \Delta$.  
  An element $\beta \in \Delta$ is maximal if and only if the singleton $\aset{\beta}$ is an open set. 
  So, the fibre $\inv{\varphi}(\aset{\beta})$ of an open $\aset{\beta}$ is a disjoint union of open sets homeomorphic to $\aset{\beta}$ if and only if every $\alpha$ in $\inv{\varphi}(\beta)$ is maximal.
  Thus, if $(\varphi, \calB)$ is a branched cover, then 
  \begin{enumerate}
    \item \label{item:FibresOfMaximal} If $\beta$ is in $\max \Delta$, then $\inv{\varphi}(\beta) \subset \max \Sigma$.
    \item We have that  $\max \Delta \subset \Delta \setminus \calB$.
  \end{enumerate}
  Moreover, if $\varphi$ fulfills \eqref{item:FibresOfMaximal}, then $(\varphi, \Delta \setminus \max \Delta)$  is a branched cover.  
\end{ex}

To extend the count of points to the fibres over $B$, we use a multiplicity on the points of $X$,
which locally agrees with an appropriate neighbourhood of the fibre.

\begin{de}[local degree] 
      \label{de:local-degree}
      Let $F: X \to Y$ be a map with finite fibres, $m_X : X \to \ZZ_{\ge1}$ an index map, and $V \subset X$. 
      We define the \emph{local degree}\index{local degree} map $\deg (F, m_X, V) : Y \to \ZZ_{\ge0}$ as
	\begin{align} \label{eq:degree-fmv}
		\deg (F, m_X, V)(y) = \sum_{\substack{ x \in \inv F(y) \\ x \in V} } m_X(x).
	\end{align}	
\end{de}

	From now on we assume that $F$ has finite fibres and $Y$ is connected.
	The convention for empty sums in Equation~\eqref{eq:degree-fmv} is to evaluate to~0. 
  The domain of $\deg (F, m_X, V)$ is $Y$ in order to make the following observation:
  for any pair of sets $V_1, V_2 \subset X$ we have that
	\begin{align} \label{eq:inclusion-exclusion}
	\deg (F, m_X, V_1) + \deg (F, m_X, V_2) = \deg (F, m_X, V_1 \cup V_2) + \deg (F, m_X,  V_1 \cap V_2).
	\end{align}
	Moreover, if $\inv F(y) \cap V_1 = \inv F(y) \cap V_2$, we have the transition equality
	\begin{align} \label{eq:transition}
    \deg (F, m_X, V_1)(y) = \deg (F, m_X, V_2)(y).
	\end{align}
	
  Finally, we say that $\deg(F, m_X, V)$ is constant if it is constant over $F(V)$.

\begin{de} 
      \label{de:IndexedBranchedCover}
  A pair $(F, m_X)$ of a branched covering $F: X \to Y$, 
  and an index map $m_X : X \to \ZZ_{\ge1}$, 
  is an \emph{indexed branched cover}\index{indexed branched cover} if 
  for every connected open set $U \subset Y$ 
  and connected component $V$ of $\inv F(U)$ 
  the local degree $\deg(F, m_X, V)$ is constant.  
\end{de}

	In particular, since $Y$ is open and connected, 
  this means that the count with multiplicity $m_X$ of the points in the fibre $\inv F(y)$ 
  is a constant $\deg (F, m_X)$ over $Y$.

\begin{ex} 
      \label{ex:LocalDegreeRealizability}
Let $f : X \to Y$ be a non-constant holomorphic map of compact Riemann surfaces, and $x$ be a point in $X$. 
Recall that there are charts for $X$ and $Y$ that express $f$ in a neighbourhood $U$ of $x$ as $f = z^{k}$.  
The integer $k$ is independent of how the charts are chosen, so the map $m_X : X \to \ZZ_{\ge1}$ sending $x \mapsto k$ is well defined. 
This is called the \emph{ramification index}\index{ramification index} of~$x$. 
The pair $(f, m_X)$ is an indexed branched cover, with a branch locus that consists of finitely-many points of $Y$. 
\end{ex}

\begin{ex} 
     \label{ex:TropicalMorphism}
	Consider the map $\Gamma \to \Delta$ from Figure \ref{fig:covers} (b). 
  The following table specifies an index map $m_\mG$ which makes $(\Gamma \to \Delta, m_\mG)$ 
  a degree-3 indexed branched cover with branch locus equal to $\aset{A,B}$.
	\begin{center}
		\begin{tabular}{c | c | c | c}
			$x$ & $m_\Gamma(x)$ & $x$ & $m_\Gamma(x)$\\
			\hline
			$A_1$ & 3 & $s_1$ & 2 \\
			$B_1$ & 3 & $s_2$ & 1 \\
			$C_1$ & 1 & $t_1$ & 1 \\
			$C_2$ & 2 & $t_2$ & 2
		\end{tabular}
	\end{center}
Note that $\aset{ s, t } \subset \aset{s,t,C}$ is indeed a dense set, since its closure is $\downset s \cup \downset t = \Delta$.
\end{ex}

\begin{re}
	Definition~\ref{de:IndexedBranchedCover} is inspired by, and synthesizes together, Definitions~III.3 and~III.4 from \cite{pay06a}.  
\end{re}

\mysubsection{A criterion for being an indexed branched cover}
      \label{sub:ACriterionForBeingAnIndexedBranchedCover}
      As a first step in the direction of Question~\ref{que:Criteria}, we show that it is enough to consider certain bases for the topology of $Y$. 
      We let $\pi_0(X)$ denote the set of connected components of $X$.

\begin{re} 
  The notation $\pi_0(X)$ should cause no confusion, since we only consider spaces that are locally path-connected, 
  i.e.~there is a basis consisting of path-connected sets, see Propositions~\ref{prop:PosetPathConnectedBasis} and~\ref{prop:CountableBasisOfPathConnectedSets}.
  In this situation the connected components and the path-connected components coincide \cite[Theorem 25.5]{mun00}.
\end{re}

\begin{lm} 
   \label{lm:LocalDegreeUnionTwoSets}
  Let $F \!:\! X \to Y$ be a branched cover, $m_X : X \to \ZZ_{\ge 1}$, and $\Aq U 1$, $\Aq U 2$ connected open sets in $Y$ with $\deg(F, m_X, \Aq V q)$ constant for $\Aq V q \in \pi_0(\inv F(\Aq U q))$, for $q = 1$ and 2. 
  The degree $\deg (F, m_X, V)$ is constant for $V \in \pi_0(\inv F(\Aq U 1 \cup \Aq U 2))$.
\end{lm}
 
\begin{proof}
We are done if $\Aq U 1$ and $\Aq U 2$ are disjoint, so assume that $\Aq U 1 \cap \Aq U 2 \ne \varnothing$.
Let $V$ be in $\pi_0(\inv F( \Aq U 1 \cup \Aq U 2 ))$. 
If $\aset{ \Aq V q_j}_{j \in J}$ is the family of connected components of $\inv F(\Aq U q)$ that intersect~$V$, for $q = 1$ and 2, we claim that 
	\[ V = \bigcup_{j \in J} \Aq V 1_j \cup \bigcup_{k \in K} \Aq V 2_k.\] 
Indeed, it is clear that the set on the right contains the set on the left, and the other containment follows from the observation that
if $S$ is a connected set with $S \cap V \ne \varnothing$, then $S \cup V$ is connected, so $S \subset V$. 
Let $y_0 \in \Aq U 1 \cap \Aq U 2$, $y_1 \in \Aq U 1$, and $y_2 \in \Aq U 2$. 
By Equations~\eqref{eq:inclusion-exclusion} and \eqref{eq:transition} we have
	\begin{align} \label{eq:Vj}
		\deg \left (F, m_X, \ \bigcup_{j \in J} \Aq V 1_j \right)(y_1) 
		 &= \sum_{j \in J} \deg \left (F, m_X, \, \Aq V 1_j \right)(y_1) \\
		\nonumber &= \sum_{j \in J} \deg \left (F, m_X, \, \Aq V 1_j \right)(y_0) \\
 		\nonumber &= \deg \left (F, m_X,  \ \bigcup_{j \in J} \Aq V 1_j \right)(y_0).
	\end{align}
	Likewise, we derive
	\begin{align} \label{eq:Vk}
		\deg \left (F, m_X, \ \bigcup_{k \in K} \Aq V 2_k \right)(y_2) = \deg \left (F, m_X, \ \bigcup_{k \in K} \Aq V 2_k \right)(y_0).  
	\end{align}

	Since $y_0$ is in  $\Aq U 1$ and $\aset{ \Aq V 1_j}_{j \in J}$ is the family of connected components of $\inv F(\Aq U 1)$ that intersects $V$, we have that 
	\[\inv F(y_0) \cap \bigcup_{j \in J} \Aq V 1_j = \inv F(y_0) \cap V. \]
	Likewise, we obtain a similar expression in relation to $\inv F(\Aq U 2)$, and conclude
	\[ \inv F(y_0) \cap \bigcup_{j \in J} \Aq V 1_j = \inv F(y_0) \cap V = \inv F(y_0) \cap \bigcup_{k \in K} \Aq V 2_k. \]
	This implies, by Equation~\eqref{eq:transition}, that the right hand sides of Equations~\eqref{eq:Vj} and \eqref{eq:Vk} are equal. Thus, so are the left hand sides, as desired.
\end{proof}

Under mild conditions for the topological spaces $X$ and $Y$, we get a countable version. 
 
\begin{lm} 
\label{lm:PathConnectedCountablyBasisSimplifies} 
  Let $F \!:\! X \to Y$ be a branched cover and $m_X : X \to \ZZ_{\ge 1}$.
  If $X$ is locally path-connected, 
  $Y$ has a countable basis $\calU$ of connected sets,
  and $\deg (F, m_X, V)$ is constant for all $V \in \pi_0(\inv F(U))$ and $U \in \calU$, 
  then $(F, m_X)$ is an indexed branched cover.
\end{lm}

\begin{proof} 
  Let $U$ be a connected open in $Y$.
  By assumption, we can write $U$ as a countable union $U = \bigcup_{q=0}^\infty U_q$ with $U_q \in \calU$.
  Suppose there is $V \in \pi_0(\inv F(U))$ such that $\deg(F, m_X, V)$ is not constant;
  i.e.~there are $y$ and~$z$ in $Y$ such that $\deg (F, m_X, V)(y) \ne  \deg(F, m_X, V)(z)$.
  Let $\aset{y_i}_{i \in I}$ and $\aset{z_j}_{j \in J}$ be $\inv F(y) \cap V $ and $\inv F(z) \cap V$, respectively, so we get 
  \begin{align} 
    \label{eq:degreenotequal}
      \deg (F, m_X, V)(y) =  \sum_{i \in I} m_X(y_i) \ne \sum_{j \in J} m_X(z_j) = \deg(F, m_X, V)(z).
  \end{align}

  Since $U$ is open and $F$ continuous, we have that $\inv F(U)$ is an open subspace of $X$.
  Being locally path-connected is inherited to open subspaces, so we have that $\inv F(U)$ is locally path-connected. 
  Thus, the connected components and the path-connected components of $\inv F(U)$ coincide, 
  so $V \in \pi_0(\inv F(U))$ is path-connected.
  So for each pair $i \in I,j \in J$ we can choose a path $P_{ij}: [0,1] \to X$ connecting $y_i$ with $z_j$.
  As $F \circ P_{ij}$ is a continuous map, the image $\im F \circ P_{ij}$ is a compact set in $Y$.
  Since $\aset{U_q}$ is a cover of $\im F \circ P_{ij}$, we can choose a finite subcover and let $r_{ij}$ be the highest index of the $U_q$ in this finite subcover.

  Let $r = \max r_{ij}$ and $U_r = \bigcup_{q=1}^r U_q$.
  Note that $\im P_{ij} \subset \inv F(U_r)$ for all pairs $i,j$.
  So the connected component $\tilde V$ of $\inv F(U_r)$ that contains  $y_1$, also contains all $y_i$, $z_j$ since $\tilde V$ is a path-connected component as well.
  Thus, by applying finitely-many times Lemma~\ref{lm:LocalDegreeUnionTwoSets}, we get that $\deg (F, m_X, \tilde V)(y) =  \deg(F, m_X, \tilde V)(z)$,
  which in particular implies
  \[ \sum_{i \in I} m_X(y_i) = \sum_{j \in J} m_X(z_j), \] 
  contradicting Equation~\eqref{eq:degreenotequal}. Thus, $\deg(F, m_X, V)$ is constant, as desired.
\end{proof}

\mysubsection{Bases of path-connected sets for \texorpdfstring{$\Sigma$}{Σ} and \texorpdfstring{$\topo \Sigma$}{|Σ|}}
    \label{sub:BasesOfPathConnectedSetsForSigmaAndTopoSigma}
Given a polyhedral complex we describe bases for its poset and topological realization that satisfy the conditions 
of Lemmas~\ref{lm:LocalDegreeUnionTwoSets} and \ref{lm:PathConnectedCountablyBasisSimplifies}.
First, for posets, recall that a \emph{principal open set} in the poset topology is a set of the form $\upset_\Sigma \alpha$. We have: 

\begin{prop}
        \label{prop:PosetPathConnectedBasis}
        Let $\Sigma$ be a finite poset. The family of principal open sets of $\Sigma$ is a finite basis of path-connected sets for the poset topology on $\Sigma$. 
\end{prop}

\begin{proof} 
  The family $\aset{\upset \alpha }_{\alpha \in \Sigma}$ is a finite basis for the poset topology. 
  To see that $\upset \alpha$ is path-connected, note that given $\gamma_1, \gamma_2$ in $\upset_\Sigma \alpha$ the following map $f: [0,1] \to \Sigma$ with $f(0) = \gamma_1$ and $f(1) = \gamma_2$ is continuous:
  \begin{align} 
        \label{eq:PathInPoset}
    f(t) =  \begin{cases} 
      \gamma_1 & \text{for } 0 \le t < \frac 1 3, \\
      \alpha & \text{for } \frac 1 3 \leq t\leq \frac 2 3, \\
      \gamma_2 & \text{for } \frac 2 3 < t \le 1. 
   \end{cases}
  \end{align}
 Thus, $\upset_\Sigma \alpha$ is path-connected, as desired.  
\end{proof}

Second, for polyhedral complexes, we glue balls in several $N^\alpha_\RR$ to obtain a topological basis of connected sets for~$\topo \Sigma$.
For this we fix an Euclidean metric $d_\alpha$ on each $N^\alpha_\RR$, and assume that all the face morphisms in $\Sigma$ are isometries.
We argue in the proof of Theorem~\ref{thm:SummarySection1} at the end of Subsection~\ref{sub:RelatingIndexedBranchedCoves} that for any polyhedral complex $\Sigma$ there is $\Sigma'$ whose face morphisms are isometries and such that $\topo \Sigma$ and $\topo {\Sigma'}$ are homeomorphic, hence this assumption can be dropped from any topological statement.

We denote by $B(N^\alpha_\RR, x, \varepsilon) = \aset{y \in N^\alpha_\RR \suchthat d_\alpha(x,y) < \varepsilon }$ the open ball in $ N^\alpha_\RR $ centred at $x$ with radius $\varepsilon$.
When $x$ is in $N^\alpha_\QQ$ and $\varepsilon$ is in $\QQ$ we say that $B(N_\RR, x, \varepsilon)$ is \emph{rational}. 
Since $\QQ$ is dense in $\RR$, the family of rational balls is a countable basis of connected sets for the Euclidean topology of $N_\RR$. 
Since all $f_{\alpha\beta}$ are isometries we have $d_\alpha = d_\beta \circ  \sigma(f_{\alpha\beta})$.
So for a point $x$ in $\topo \Sigma$ we define the \emph{principal ball centred at $x$ with radius $\varepsilon$} as
\begin{align}
  \label{eq:PrincipalBall}
  \upsetballsub \Sigma x \varepsilon = \bigcup_{\beta \succeq \poly_\Sigma(x)} p_\beta(B(N^\beta_\RR, \inv p_\beta(x), \varepsilon) \cap \sigma_\beta). 
\end{align}

For an open set $U \subset \topo \Sigma$ with $x \in U$ and small enough $\varepsilon$, the set $\upsetballsub \Sigma x \varepsilon$ is an open neighbourhood.
This follows from two lemmas. 
  
\begin{lm} 
  \label{lm:preimagePrincipalBall}
  Let $\Sigma$ be a polyhedral complex whose face morphisms are isometries.
  For every  $x$ in $\topo \Sigma$ and $\gamma$ in $\upset_{\Sigma} \poly_\Sigma x$
  we have that 
  \begin{align}
    \label{eq:SmallDistance}
  \inv p_\gamma \left(    \upsetballsub \Sigma x \varepsilon   \right) 
    =  B(N^\gamma_\RR, \inv p_\gamma(x), \varepsilon) \cap \sigma_\gamma. 
  \end{align}
\end{lm}

\begin{proof} 
  The set on the left contains the right one because of the assumption that $\gamma \succeq \poly_\Sigma (x)$.
For the other containment, let $\hat y$ be in $\inv p_\gamma(\upsetballsub \Sigma x \varepsilon)$, namely for some $\beta \in \upset \poly_\Sigma(x)$ we have that
\[   p_\gamma(\hat y) \in p_\beta(B(N^\beta_\RR, \inv p_\beta(x), \varepsilon) \cap \sigma_\beta). \]
Let $y$ equal $p_\gamma(\hat y)$, and note that this point is in $\im p_\gamma \cap p_\beta$, so by Lemma~\ref{lm:impalphaIntersectsimpbeta} there exists $\eta$ in $\Sigma$ with morphisms $f: \eta \to \beta$ and $g : \eta \to \gamma$, such that $y$ is in $\im p_\eta$.
As $f, g$ are isometries and $f = \inv p_\beta \circ p_\eta$, $g = \inv p_\gamma \circ p_\eta$, we calculate
\begin{align*} 
  d_\gamma(\hat y, \inv p_\gamma(x)) &= d_\gamma(\inv p_\gamma(y), \inv p_\gamma(x)) = d_\gamma(\inv p_\gamma \circ p_\eta \circ \inv p_\eta(y), \inv p_\gamma \circ p_\eta \circ \inv p_\eta(x)) \\
  &= d_\gamma(g \circ \inv p_\eta(y), g \circ \inv p_\eta(x)) = d_\eta(\inv p_\eta(y), \inv p_\eta(x)) \\
  & = d_\beta( \inv p_\beta(y), \inv p_\beta(x)) \le \varepsilon, 
\end{align*}
and conclude that $\hat y$ is in  $B(N^\gamma_\RR, \inv p_\gamma(x), \varepsilon)$, as desired. 
\end{proof}

Second, we express containment in terms of the universal maps $p_\alpha$.
Denote by $\closure -$ the closure of a set. 

\begin{lm} 
   \label{lm:containmentAndUniversalMaps} 
   Let $\Sigma$ be a polyhedral complex, $U$ and $V$ subsets of $\topo \Sigma$, and $\calU \subset \Sigma$ such that $\poly_\Sigma U \subset \closure \calU$. 
   If $\inv p_\beta(U) \subset \inv p_\beta(V)$ for all $\beta$ in $\calU$, then $U \subset V$. 
\end{lm}

\begin{proof} 
  Let $x$ be in $U$, set $\alpha = \poly_\Sigma x$, and $\hat x = \inv p_\alpha(x)$. 
  Since $\poly_\Sigma U \subset \closure \calU$, 
  there is a $\beta$ in $\calU$ and a morphism $f : \alpha \to \beta$ in $\Sigma$.
  Hence, $x = p_\alpha(\hat x) = p_\beta \circ f(\hat x)$ gives that $f(\hat x) \in \inv p_\beta(U)$, so by assumption also in $\inv p_\beta(V)$, hence $x$ is in $V$ as desired. 
\end{proof}

Putting both ingredients together, we get:

\begin{lm} 
  \label{lm:smallEpsilonGivesContainment}
  Let $\Sigma$ be a polyhedral complex whose face morphisms are isometries, 
  $V \subset \topo \Sigma$ an open set, and $x$ in~$V$.
  There is a constant $K$ such that if $\varepsilon < K$ we have
  \[ \upsetballsub \Sigma x \varepsilon \subset V.  \]
\end{lm}

\begin{proof} 
  Let $\alpha = \poly_\Sigma x$, set $U = \upsetballsub \Sigma x \varepsilon $ and $\calU = \upset_\Sigma \alpha$.
  By Equation~\ref{eq:PrincipalBall} we have that $\poly_\Sigma U \subset \closure \calU$. 
  Since $V$ is open, for every $\beta \in \calU$ we have that $\inv p_\beta(V)$ is open.
  Moreover, note that $\inv p_\beta(x) \in \inv p_\beta(V)$.
  Thus, we can choose $\varepsilon_\beta$ such that 
  \[B(N^\beta_\RR, \inv p_\beta(x), \varepsilon_\beta) \cap \sigma_\beta \subset \inv p_\beta(V).\]
  By Lemma~\ref{lm:preimagePrincipalBall} the left hand side equals $\inv p_\beta(\upsetballsub \Sigma x {\varepsilon_\beta})$.
  Hence, if we take $K = \min_{\beta \in \calU} \varepsilon_\beta$, we are done by Lemma~\ref{lm:containmentAndUniversalMaps}.
\end{proof}


\begin{lm} 
  \label{lm:smallEpsilonImpliesOpen}
  Let $\Sigma$ be a polyhedral complex whose face morphisms are isometries, 
  and $x$ in~$\topo \Sigma$.
  There is a constant $K$ such that if $\varepsilon < K$ we have that $\upsetballsub \Sigma x \varepsilon$ is open.
\end{lm}

\begin{proof} 
  Let $\alpha = \poly_\Sigma$.
  Applying Lemma~\ref{lm:smallEpsilonGivesContainment} to $\inv \poly_\Sigma(\upset_\Sigma \alpha)$, 
  which is open because $\poly_\Sigma$ is continuous,  
  we get a $K$ such that if $\varepsilon < K$, then $\upsetballsub \Sigma x \varepsilon \subset \inv \poly_\Sigma(\upset_\Sigma \alpha)$.
  So $\poly_\Sigma (\upsetballsub \Sigma x \varepsilon)= \upset_\Sigma \alpha$, and we are done by Lemma~\ref{lm:preimagePrincipalBall}. 
\end{proof}

Note that any point $y$ in $\upsetballsub \Sigma x \varepsilon$ is connected by a path to $x$. 
Also, the set of rational points  
$\Sigma_\QQ = \bigcup_{\alpha \in \Sigma} p_\alpha(N^\alpha \otimes_\ZZ \QQ \cap \sigma_\alpha)$ 
of $\Sigma$ is countable and dense. So we arrive to:

\begin{prop}
  \label{prop:CountableBasisOfPathConnectedSets}
  Let $\Sigma$ be a polyhedral complex whose face morphisms are isometries. 
  The family of principal balls of $\Sigma$ that are rational and open, i.e.~$\upsetballsub \Sigma x \varepsilon$ with $x \in \Sigma_\QQ$ and $\varepsilon \in \QQ$ small enough, 
  is a countable basis of path-connected sets for $\topo \Sigma$. \qed
\end{prop}

\mysubsection{\Combinatorial{} morphisms}
    \label{sub:CombinatorialMorphisms}

    Let $\Phi : [\sigma : \Sigma \to \polyintface] \to [\delta : \Delta \to \polyintface]$ be a morphism in $\RPC$, and $\topo \Phi : \topo \Sigma \to \topo \Delta$ its topological realization.
We give a condition which implies that $\topo \Sigma$ arises as several copies of $\topo \Delta$ glued together in a manner prescribed solely by $\varphi : \Sigma \to \Delta$.
This strengthens Equation~\eqref{eq:containment-fibres} to a bijection, 
which is crucial to relate the count of points in the fibres of $\topo \Phi$ and $\varphi$.

\begin{de}
	We say that $\Phi = (\varphi, \aset{\Phi_\alpha }_{\alpha \in \Sigma} )$ is \emph{\combinatorial{} } if for all $\alpha \in \Sigma$ 
  we have that 
  $\Phi_\alpha(\sigma_\alpha) = \delta_{\varphi(\alpha)}$
  and
  $\dim \sigma_\alpha = \dim \delta_{\varphi(\alpha)}$.
\end{de}

The condition $\Phi_\alpha(\sigma_\alpha) = \delta_{\varphi(\alpha)}$ implies that $\Phi_\alpha : N^\alpha_\RR \to N^{\varphi(\alpha)}_\RR$ is a surjective linear map. 
Combined with the dimension condition, we get that the linear map $\Phi_\alpha$ is bijective, hence a homeomorphism. 
In fact, it is straightforward to see that $\Phi$ is \combinatorial{} if and only if $\Phi_\alpha$ is injective and $\Phi_\alpha(\sigma_\alpha) = \delta_{\varphi(\alpha)}$ for all $\alpha$ in $\Sigma$.

\begin{re}
	In \cite{pay06a} 
  a \combinatorial{} morphism of cone complexes satisfies that $\Phi_\alpha$ maps $(N^\alpha, \sigma_\alpha)$ isomorphically to $( N^{\varphi(\alpha)}, \delta_{\varphi(\alpha)})$. 
	This is equivalent to $\Phi_\alpha :  N^\alpha_\RR \to  N^{\varphi(\alpha)}_\RR$ being injective, 
      $\Phi_\alpha(\sigma_\alpha) = \delta_{\varphi(\alpha)}$, 
      and $\Phi_\alpha(N^\alpha) = N^{\varphi(\alpha)}$. 
  We have omitted the last condition, since the index $[N_{\varphi(\alpha)} : \Phi_\alpha(N_\alpha)]$ encodes relevant information.
\end{re}

We give a combinatorial characterization in terms of the map $\varphi$.

\begin{lm} \label{lm:characterization-combinatorial-morphism}
	A morphism $(\varphi : \Sigma \to \Delta, \aset{\Phi_\alpha }_{\alpha \in \Sigma} )$ of polyhedral complexes is \combinatorial{} if and only if $\varphi$ maps $\downset_\Sigma \alpha$ isomorphically to $\downset_\Delta \varphi(\alpha)$ for all $\alpha \in \Sigma$.
\end{lm}

\begin{proof}
	Suppose that $(\varphi, \aset{\Phi_\alpha }_{\alpha \in \Sigma} )$ is \combinatorial{}, so $\Phi_\alpha :  N^\alpha_\RR \to  N^{\varphi(\alpha)}_\RR$ is 
  a bijective linear map, hence the polyhedral structure of $\sigma_\alpha$ is isomorphic to that of $\Phi_\alpha(\sigma_\alpha)$.
  By supposition $\Phi_\alpha(\sigma_\alpha) = \delta_{\varphi(\alpha)}$, so their polyhedral structure is also isomorphic. 

	For the converse, suppose that $\varphi$ maps $\downset_\Sigma \alpha$ isomorphically to $\downset_\Delta \varphi(\alpha)$ for all $\alpha \in \Sigma$. 
  We must show that   $\Phi_\alpha(\sigma_\alpha) = \delta_{\varphi(\alpha)}$
  and
  $\dim \sigma_\alpha = \dim \delta_{\varphi(\alpha)}$ for all $\alpha \in \Sigma$.
  By definition we already have $\Phi_\alpha(\sigma_\alpha) \subset \delta_{\varphi(\alpha)}$.
  For the other containment, since $\varphi(\downset_\Sigma \alpha) = \downset_\Delta \varphi(\alpha)$, we have that all the 1-faces of $\delta_{\varphi(\alpha)}$ are in~$\Phi_\alpha(\sigma_\alpha)$. 
  We are done because $\delta_{\varphi(\alpha)}$ is the convex hull of its 1-faces.
  Moreover, the equality on dimension is true because $\dim \sigma_\alpha$ equals the number of elements in a maximal chain of $\downset_\Sigma \alpha$; 
  similarly for  $\downset_\Delta \varphi(\alpha)$ and $\dim \delta_{\varphi(\alpha)}$;
  since $\downset_\Sigma \alpha$ is isomorphic to $\downset_\Delta \varphi(\alpha)$, all four numbers are the same.
\end{proof}

\begin{de} 
      \label{de:CombinatorialMorphism}
    We call a morphism of posets $\varphi : \Sigma \to \Delta$ \emph{\combinatorial{} } if and only if 
    $\Sigma$ and $\Delta$ are finite and $\varphi$ maps  $\downset_\Sigma \alpha$ isomorphically to $\downset_\Delta \varphi(\alpha)$ for all $\alpha \in \Sigma$.
\end{de}

So Lemma~\ref{lm:characterization-combinatorial-morphism} says that a morphism 
$\Phi = (\varphi : \Sigma \to \Delta, \aset{\Phi_\alpha }_{\alpha \in \Sigma} )$ in $\RPC$ is \combinatorial{} 
if and only if the underlying morphism of categories $\varphi : \Sigma \to \Delta$ is \combinatorial{}.

\subsection{Fibre products}
    \label{sub:FibreProducts}
Now we prove that when $\dtmor$ is \combinatorial{}, the space $\topo \Sigma$ is homeomorphic to the fibre product $\topo \Delta \times_\Delta \Sigma$ 
under the maps $\poly_\Delta$ and $\varphi$. 
Recall that given 
$f \colon X \to Z$
and
$g \colon Y \to Z$
in $\Top$,
the fibre product $X \times_Z Y$ is the limit of the diagram $X \rightarrow Z \leftarrow Y$,
homeomorphic to the subspace 
\[
    X \times_Z Y = \aset{ (x, y) \suchthat f(x) = g(y) }
\]
of the product topology $X \times Y$. 
The universal property states that for any $W$ with maps $p \colon W  \to X$ and $q \colon W \to Y$, there is a unique map $\phi : W \to X \times_Z Y$ that gives factorizations $p = \pi_X \circ \phi$ and $q = \pi_Y \circ \phi$, with $\pi_X$ and $\pi_Y$ the canonical projections. 

\begin{lm} \label{lm:combinatorialImpliesFibredProduct}
If $\Phi : \Sigma \to \Delta$ in $\RPC$ is \combinatorial{}, then the map 
\begin{align*}
  \phi \colon \topo \Sigma  &\to \topo \Delta \times \Sigma,  \\
  x &\mapsto (\topo \Phi(x),\ \poly_{\Sigma}(x)),
\end{align*}
is a bijection onto $\topo \Delta \times_\Delta \Sigma$.
\end{lm}

\begin{proof}
  That the image of $\phi$ is contained in $ \topo \Delta \times_\Delta \Sigma$  follows from Lemma~\ref{lm:superDiagram}.
  Let $(y, \gamma)$ be such that $\poly_{\Delta}(y) = \varphi(\gamma)$, so $y \in p_{\varphi(\gamma)}(\delta_{\varphi(\gamma)}^{\circ})$. 
  Since $(\varphi, \aset{\Phi_{\alpha}}_{\alpha \in \Sigma})$ is \combinatorial{}, the map $\Phi_{\gamma}: \topo{\sigma_{\gamma}} \to \topo{\delta_{\varphi(\gamma)}}$ is a homeomorphism. Thus, there is exactly one $x$ in $p_{\gamma}(\sigma_{\gamma}^{\circ}) = \inv \poly_{\Sigma}(\gamma)$ such that $\topo \Phi(x)=y$. This shows that $\phi$ is bijective.
\end{proof}

Given $y$ in $\topo \Delta$, the fibre $\inv{\pi_{\topo \Delta}}(y)$ is equal to $\aset{(y,\, \alpha)  \suchthat  \varphi(\alpha) = \poly_\Delta(y)}$.
Thus, Lemma~\ref{lm:combinatorialImpliesFibredProduct} strengthens Equation~\eqref{eq:containment-fibres} to an equality:

\begin{lm} 
  \label{lm:CombinatorialImpliesFibreBijection}
    If $(\varphi, \aset{\Phi_{\alpha}}_{\alpha \in \Sigma}) : \Sigma \to \Delta$ in $\RPC$ is \combinatorial{}, then for any $y$ in $\topo \Delta$ the map $\poly_{\Sigma}$ induces a bijection 
    \[ 
    \pushQED{\qed} 
        \inv{\topo \Phi}(y) \to \inv{\varphi}(\poly_{\Delta}(y)). \qedhere
    \popQED
  \]
\end{lm}

Lemma~\ref{lm:combinatorialImpliesFibredProduct} implies that $\phi$ restricts to a bijective continuous map from $\topo \Sigma$ to $\topo \Delta \times_\Delta \Sigma$.
It is left to prove that $\phi$ is an open map.
Given an open set $V$ of  $\topo \Sigma$, we give an open set $U$ of $\topo \Delta \times \Sigma$ such that $\phi(V) = U \cap (\topo \Delta \times_\Delta \Sigma)$.

Since $\poly_{\Sigma}$ is an open map, $\poly_{\Sigma}(V)$ is open. 
Also, the maps $\aset{\Phi_{\alpha}}_{\alpha \in \Sigma}$ which glue together to form $\topo \Phi$ are open. 
Thus, a natural candidate for $U$ is $\topo \Phi(V) \times \poly_\Sigma(V)$.

\begin{ex} 
  \label{ex:TooManyComponents}
  Let $(\varphi, \aset{\Phi_\alpha}_{\alpha \in \Sigma})$ be a \combinatorial{} morphism in $\RPC$ mapping a complex with three cones down to a complex with two cones. 
  Figure~\ref{fig:TooManyComponents} displays the topological realization $\topo \Phi$ and the order preserving map $\varphi$. 
  Consider an open $W \subset \topo \Delta$ 
  and a connected component $V$ of $\inv{\topo \Phi}(W)$ in the cone~$\beta_2$, as in the figure.
  We have that $U = \topo \Phi(V) \times \poly_\Sigma(V) = \topo \Phi(V) \times \aset{B_2, \beta_2 }$, but clearly $\topo \Phi(V)$ is not open.
  Indeed, since $\poly_\Delta$ is an open map, if $\topo \Phi(V)$ were open, 
  then $\poly_\Delta \circ \topo \Phi(V)$ would be open, 
  but this equals $\dtmor \circ \poly_\Sigma(V) = \aset{B, \beta}$ which is not an up-set.
  Here the proposed $U$ fails to be an open set. 
  Also, in the context of Proposition~\ref{myprop:IndexedBranchedCoversAreOpenMaps}, note that $\dtmor$ is not an open map.
  \begin{figure} \center
	\begin{minipage}{0.4\textwidth}
		\scalebox{0.8}{
		\begin{overpic}{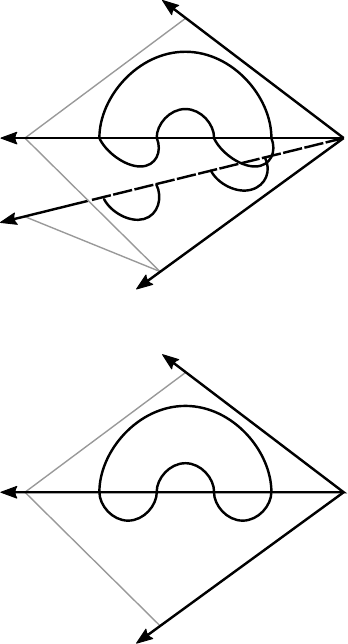} 
			\put (67,53) {\scalebox{1.3}{$\downarrow$}}
			
			\put (65,77) {\scalebox{1}{$\topo \Sigma$}} 
			\put (65,22) {\scalebox{1}{$\topo \Delta$}}

			\put (10,89) {\scalebox{0.8}{$\alpha_1$}}
			\put (4,71) {\scalebox{0.8}{$\beta_1$}} 
			\put (7,59) {\scalebox{0.8}{$\beta_2$}} 

			\put (13,36) {\scalebox{0.8}{$\alpha$}} 
			\put (10,10) {\scalebox{0.8}{$\beta$}} 

    	\put (19,100) {\scalebox{0.9}{$A_1$}}
    	\put (25,64) {\scalebox{0.9}{$V$}}
			\put (-7,78) {\scalebox{0.9}{$B_1$}} 
			\put (16,53) {\scalebox{0.9}{$C_1$}} 
			\put (-7,63) {\scalebox{0.9}{$B_2$}} 
			\put (55,77) {\scalebox{0.9}{$O_1$}} 

			\put (21,44) {\scalebox{0.9}{$A$}} 	
			\put (-5,22) {\scalebox{0.9}{$B$}} 
      \put (16,-2) {\scalebox{0.9}{$C$}} 
			\put (55,22) {\scalebox{0.9}{$O$}} 
			\put (28,30.5) {\scalebox{0.9}{$W$}} 	
		\end{overpic}
	}
	\end{minipage} 
\begin{minipage}{0.4\textwidth} \center
	\scalebox{0.8}{
		 \begin{tikzpicture}
		\node (A1) at (0,0) {$A_1$};
		\node (B1) at (1.6,0) {$B_1$};
		\node (B2) at (2.4,0) {$B_2$};
		\node (C1) at (4.2,0) {$C_1$};
    \node (O1) at (2,-1) {$O_1$};
		\node (t1) at (3.5,1.2) {$\beta_2$};
		\node (t2) at (2.5,1.2) {$\beta_1$};
		\node (s1) at (0.6,1.2) {$\alpha_1$};
    \draw (B2) -- (O1) -- (A1) -- (s1) -- (B1);
    \draw (B2) -- (t1) -- (C1);  
		\draw (B1) -- (O1) -- (C1) -- (t2);  
		 \draw[preaction={draw=white, -,line width=6pt}] (B1) -- (t2);
		 \draw[preaction={draw=white, -,line width=6pt}] (t1) -- (B2);
		 \node (A) at (0,-3.1) {$A$};
		 \node (B) at (2,-3.1) {$B$};
		 \node (C) at (4,-3.1) {$C$};
		 \node (t) at (3,-2.1) {$\beta$};
		 \node (s) at (1,-2.1) {$\alpha$};
     \node (O) at (2,-4.1) {$O$};
     \draw (O) -- (A) -- (s) -- (B) -- (t) -- (C) -- (O);
     \draw (O) -- (B);
		 \node (t2) at (5,0) {\scalebox{1.25}{$\Sigma$}};
		 \node (t2) at (5,-1.5) {\scalebox{1.4}{$\downarrow$}};
		 \node (t) at (5,-2.9) {\scalebox{1.25}{$\Delta$}};
		 \end{tikzpicture}
	}
\end{minipage}
\caption{\label{fig:TooManyComponents} On the left, the topological realization $\topo \Phi$ of a \combinatorial{} morphism in $\RPC$, an open subset $W$ of the codomain $\topo \Delta$, and the fibre $\inv{\topo \Phi}(W)$ above $W$. On the right the map $\varphi : \Sigma \to \Delta$ of face posets. } 
\end{figure}
\end{ex}


So the proof that $\phi$ is open, which implies that $\phi$ is a homeomorphism, is more involved.


\begin{lm} 
    \label{lm:phiIsOpen}
    Let $\Phi : \Sigma \to \Delta$ in $\RPC$ be \combinatorial{}.
    Assume that $\Sigma$ and $\Delta$ are glued by face morphisms that are isometries.
    The map $\phi$ given by $x \mapsto (\topo \Phi(x), \poly_\Sigma (x))$ is an open map.
\end{lm}

\begin{proof} 
  Let $V$ be an open subset of $\topo \Sigma$, 
  choose any point $x$ in $V$, 
  set $y = \topo \Phi(x)$, 
  $\alpha = \poly_\Sigma(x)$ 
  and $\beta = \varphi(\alpha) = \poly_\Delta(y)$.
  We are done if we exhibit an open neighbourhood of $\phi(x) = (y, \alpha)$ 
  contained in $\phi(V)$.
  Note that by Lemma~\ref{lm:combinatorialImpliesFibredProduct} we have
  \begin{align} 
    \label{eq:phiV}
    \phi(V) = (\topo \Phi(V) \times \poly_\Sigma(V)) \cap (\topo \Delta \times_\Delta \Sigma ). 
  \end{align}
  We claim that for small enough $\varepsilon$, the neighbourhood of $(y,\alpha)$ given by
  \[ U_\varepsilon  = (\upsetballsub \Delta y {\varepsilon} \times \upset_\Sigma \alpha) \cap (\topo \Delta \times_\Delta \Sigma), \]
  is an open and contained in $\phi(V)$.

  By Proposition~\ref{prop:CountableBasisOfPathConnectedSets} there is $\varepsilon_x > 0$ such that $\upsetballsub \Sigma x {\varepsilon_x}$ is open and contained in $V$;
  there is also $K_x$ such that if $0 < \varepsilon < K_x$ the set $\upsetballsub \Delta y \varepsilon$ is open.
  Let $(z,\eta)$ be in $U_\varepsilon$.
  Since $(z,\eta)$ is in $\topo \Delta \times_\Delta \Sigma$, we have that $\gamma = \varphi(\eta)$ equals $\poly_\Delta(z)$.  
  We also have that $\eta$ is in $\upset_\Sigma \alpha$, so Lemma~\ref{lm:preimagePrincipalBall} gives
  \begin{align*} 
    \inv p_\eta( \upsetballsub \Sigma x {\varepsilon_x} ) = B(N^\eta_\RR, \inv p_\eta(x), \varepsilon_x) \cap \sigma_\eta.
  \end{align*}
  Since $\Phi$ is \combinatorial{}, $\Phi_\eta$ maps $\sigma_\eta$ homeomorphically to $\delta_\gamma$.
  As $  B(N^\eta_\RR, \inv p_\eta(x), \varepsilon_x) \cap \sigma_\eta $ is open in $\sigma_\eta$, 
  the set $\Phi_\eta(   B(N^\eta_\RR, \inv p_\eta(x), \varepsilon_x) \cap \sigma_\eta  )$ is open in $\delta_\gamma$,
  so there is $\varepsilon_\eta$ such that 
  \begin{align*} 
    B(N^\gamma_\RR, \Phi_\eta \circ \inv p_\eta(x), \varepsilon_\eta) \cap \delta_\gamma \subset  \Phi_\eta(   B(N^\eta_\RR, \inv p_\eta(x), \varepsilon_x) \cap \sigma_\eta  ) \subset \inv p_\gamma( \topo \Phi(  \upsetballsub \Sigma x {\varepsilon_x} ) ). 
  \end{align*}
  The second containment above follows from Diagram~\ref{dia:genrpc-morphism};
  from the same diagram we also get that $  \Phi_\eta \circ \inv p_\eta(x) = \inv p_\gamma(y)$. So if $\varepsilon < \varepsilon_\eta$, we get
  \begin{align} 
    \label{eq:zIsInTheSet} 
    z \in p_\gamma( B(N^\gamma_\RR, \inv p_\gamma(y), \varepsilon) \cap \delta_\gamma) \subset \topo \Phi(  \upsetballsub \Sigma x {\varepsilon_x}). 
  \end{align}

  If we choose an $\varepsilon_\eta$ for each $\eta$ in $\upset_\Sigma \alpha$ as above, 
  and set  $\varepsilon$ to be the minimum of the $\varepsilon_\eta$ and $K_x$, 
  we have by Equation~\eqref{eq:zIsInTheSet} and Equation~\eqref{eq:PrincipalBall}
  that 
  \[(\upsetballsub \Delta y \varepsilon \times \upset_\Sigma \alpha) \cap (\topo \Delta \times_\Delta \Sigma) \subset (\topo \Phi(\upsetballsub \Sigma x {\varepsilon_x}) \times \upset_\Sigma \alpha ) \cap (\topo \Delta \times_\Delta \Sigma).  \]
  Note that $\topo \Phi( \upsetballsub \Sigma x {\varepsilon_x})$ is contained in $\topo \Phi(V)$. 
  Also, since $V$ is open,  $\poly_\Sigma(V)$ is an up-set, so $\upset_\Sigma \alpha$ is contained in it.
  Comparing with Equation~\eqref{eq:phiV} we see that $U_\varepsilon \subset \phi(V)$; and as both $\upsetballsub \Delta y \varepsilon$ and $\upset_\Sigma \alpha$ are open, we are done. \qedhere
\end{proof}

\mysubsection{Relating indexed branched covers}
    \label{sub:RelatingIndexedBranchedCoves}
Let $\Sigma$ and $\Delta$ be polyhedral complexes, $\Phi : \Sigma \to \Delta$ a \combinatorial{} morphism and $\topo \Phi$ its topological realization.
We look at connected components of preimages in $\Sigma$ of principal open sets $\upset_\Delta \beta \subset \Delta$,
and of preimages in $\topo \Sigma$ of open principal balls in $\topo \Delta$. 
With Lemma~\ref{lm:PathConnectedCountablyBasisSimplifies} this culminates the proof of Theorem~\ref{thm:SummarySection1}.

\begin{lm}
  \label{lm:connected-components-combinatorial-morphism}
If $\varphi : \Sigma \to \Delta$ is a \combinatorial{} morphism of posets and $\beta$ is in $\Delta$, then 
\[ \conncomp {\inv{\varphi}(\upset_\Delta \beta )} = \aset{\upset_\Sigma \alpha }_{\alpha \in \inv{\varphi}(\beta)}. \]
\end{lm}

\begin{proof}
  First observe that for any $\gamma \in \inv{\varphi}(\upset_\Delta \beta  )$ we have that~$\beta \in \downset_\Delta \varphi(\gamma)$. 
  Since $\varphi$ is a \combinatorial{} morphism, it maps $\downset_\Sigma \gamma$ bijectively to $\downset_\Delta \varphi(\gamma )$. 
  Thus, there is an $\alpha \in \Sigma$ such that $\alpha \preceq \gamma$ and $\varphi(\alpha) = \beta$,
  which gives	
  \begin{align} \label{eq:ConnectedComponentsOfFibre}
	\inv \varphi(\upset_\Delta \beta  ) = \bigcup_{\alpha \in \inv \varphi(\beta)} \upset_\Sigma \alpha.
	\end{align}
  This union is disjoint, because if there were $\alpha, \alpha'$ in $\inv \varphi(\beta)$ such that $\upset_\Sigma \alpha \cap  \upset_\Sigma \alpha'$ contained an element $\gamma$, 
	then $\varphi$ would map $\downset_\Sigma \gamma$ isomorphically to $ \downset_\Delta \varphi(\gamma )$.  
  Since $\alpha, \alpha' \in \downset_\Sigma \gamma$, this gives $\varphi (\alpha) \ne \varphi (\alpha')$, a contradiction. 
  Finally, any open set that contains $\alpha$ must contain $\upset_\Sigma \alpha$, so each $\upset_\Sigma \alpha$ in Equation~\eqref{eq:ConnectedComponentsOfFibre} is connected.
\end{proof}


The following is a version of Lemma~\ref{lm:CombinatorialImpliesFibreBijection} for connected components.

\begin{prop} 
  \label{prop:FibreOfPrincipalBall}
  Let $\Phi$ be \combinatorial{}, $\Sigma$ and $\Delta$ have face morphisms that are isometries,
  $U = \upsetballsub \Delta y \varepsilon \subset \topo \Delta$ be an open principal ball, and $\calU = \poly_\Delta U$.
  The map $\poly_\Sigma$ induces a bijection of connected components of the fibres
  \[ 
	\conncomp{\inv {\topo \Phi} (U)} 
	\to 
	\conncomp {\inv \varphi( \calU )}.
	\]
\end{prop}

\begin{proof} 
  Let $V$ be a component in~$\conncomp{\inv {\topo \Phi}(U)}$.
  We first show that $\poly_\Sigma(V)$ is a component in $\conncomp {\inv \varphi( \calU )}$.
  Set $\beta = \poly_\Delta y$. 
  Since $U$ is open we get $\calU = \poly_\Delta U = \upset_\Delta \beta$.
  Assume there is a point $x$ in $V$ such that $\topo \Phi (x) = y$ and set $\alpha = \poly_\Sigma(x)$.
  Since $\poly_\Sigma$ is a continuous open map, we have that $\poly_\Sigma(V)$ is a connected open set containing $\upset_\Sigma \alpha$. 
  Since $\alpha$ is in  $\inv \varphi(\beta)$,
  Lemma~\ref{lm:connected-components-combinatorial-morphism} implies that $\upset_\Sigma \alpha$ is a connected component of $\conncomp {\inv \varphi( \calU )} = \conncomp {\inv \varphi( \upset_\Delta \beta )}$,  so $\poly_\Sigma(V) = \upset_\Sigma \alpha$.

  To prove there is a point $x$ in $V$ such that $\topo \Phi (x) = y$, 
  choose any $z$ in $V$ and set $\gamma = \poly_\Sigma(z)$.
  We have that $\beta \preceq \varphi(\gamma)$, beause $\topo \Phi(V) \subset U$ and $\varphi \circ \poly_\Sigma = \poly_\Delta \circ \topo \Phi$.
  So there is a morphism $f_{\beta \varphi(\gamma)} : \delta_\beta \to \delta_{\varphi(\gamma)}$, and we have $p_\beta = p_{\varphi(\gamma)} \circ f_{\beta \varphi(\gamma)}$.
  Since $y$ is in $\im p_\beta$, we have $y$ in $\im p_{\varphi(\gamma)}$, and we may set $\hat y = \inv p_{\varphi(\gamma)}(y)$.

  As $\Phi$ is \combinatorial{}, the affine map $\Phi_\gamma$ sends $\sigma_\gamma$ homeomorphically to $\delta_{\varphi(\gamma)}$.
  So $\hat y$ is in $\im \Phi_\gamma$, and we may set $\hat x = \inv \Phi_\gamma (\hat y)$.
  By Diagram~\ref{dia:HowTopologicalRealizationGlues}, we have that $\topo \Phi \circ p_\gamma(\hat x) = p_{\varphi(\gamma)} \circ \Phi_{\gamma}(\hat x) = p_{\varphi(\gamma)}(\hat y) = y$.
  Set $x = p_\gamma (\hat x)$, so $\topo \Phi(x) = y$. 
  Note that $x$ is in $p_\gamma(\inv \Phi_\gamma(  \inv p_{\varphi(\gamma)} (U) ))$.
  We claim the latter set is contained in $V$.

  Since $\poly_\Delta y = \beta \preceq \varphi(\gamma)$, we may apply Lemma~\ref{lm:preimagePrincipalBall} to get that
  \[
    \inv p_{\varphi(\gamma)} (U) = B(N^{\varphi(\gamma)}_\RR, \inv p_{\varphi(\gamma)}(y), \varepsilon) \cap \delta_{\varphi(\gamma)},
  \]
  which in particular means that $ \inv p_{\varphi(\gamma)} (U)$  is a convex set, hence connected.
  As $\Phi_\gamma$ is a homeomorphism,  $\inv \Phi_\gamma(  \inv p_{\varphi(\gamma)} (U) )$ is connected as well, and so is $p_\gamma( \Phi_\gamma(  \inv p_{\varphi(\gamma)} (U) ))$.
  Finally, recall we have chosen a point $z$ in $V$ and set $\gamma = \poly_\Sigma(z)$.
  Since $\topo \Phi(z) \in \topo \Phi(V) \subset U$ and $\poly_\Sigma(z) = \gamma$, by Diagram~\ref{dia:HowTopologicalRealizationGlues}  we have that $z$ is in $p_\gamma(\inv \Phi_\gamma(  \inv p_{\varphi(\gamma)} (U) ))$.
  So $p_\gamma(\inv \Phi_\gamma(  \inv p_{\varphi(\gamma)} (U) )) \subset \inv {\topo \Phi}(U)$ is connected, and intersects the connected component $V$, so we conclude the desired inclusion.
  
  To conclude, we have proven that each component of $\inv {\topo \Phi}(U)$ intersects $\inv {\topo \Phi}(y)$,  for $y$ the centre of the ball $U$.
  Moreover, we have proven that for $x$ in $\inv {\topo \Phi}(y)$ and $V$ the component containing $x$, we get $\poly_\Sigma(V) = \upset_\Sigma( \poly_\Sigma(x))$.
  These two facts, together with Lemma~\ref{lm:CombinatorialImpliesFibreBijection} and Lemma~\ref{lm:connected-components-combinatorial-morphism} imply surjectivity and injectivity. 
\end{proof}


  \begin{re} 
    \label{re:RestrictionsOnUNecessary}
    In Figure~\ref{fig:TooManyComponents} the set $U$ is connected, 
    and the set  $\calU = \poly_\Delta U$ is principal,   
    yet $\conncomp{\inv {\topo \Phi}(U) }$ has three elements, while $\conncomp{\inv \varphi (\calU)}$ has only two.
    This highlights the need to impose restrictions on $U$ in Proposition~\ref{prop:FibreOfPrincipalBall}. 
  \end{re}
   
Now, we show that if one of the vertical maps in Diagram~\ref{dia:super-diagram} is an indexed branched cover, then so is the other. 
Consider an index map $m_\Phi: \topo \Sigma \to \ZZ_{\ge 1}$,
a connected open~$U \subset \topo \Delta$, 
a connected component $V$ in $\conncomp { \inv {\topo \Phi}(U) }$, 
and a point $y$ in $U$. 
We have that
\[ 	\deg( \topo \Phi, m_\Phi, V)(y) 
= \sum_{\substack{ x \in \inv {\topo \Phi}(y) \\ x \in V} } m_\Phi(x).\]
Likewise, for an index map $m_\varphi : \Sigma \to \ZZ_{\ge 1} $, a connected open $\upset_\Delta \beta \subset \Delta$, a connected component $\calV$ in $\conncomp { \inv {\topo \Phi}(\upset_\Delta \beta) }$, and a point $\gamma$ in $\upset_\Delta \beta$.
We have that
\[ 	\deg( \varphi, m_\varphi, \calV)(\gamma) 
= \sum_{\substack{ \alpha \in \inv \varphi(\gamma) \\ \alpha \in \calV} } m_\varphi(\alpha).\]
Our aim is to relate both degrees.
The first step compares the index sets $\inv {\topo \Phi}(y) \cap V$ and $\inv \varphi(\gamma) \cap \calV$ by putting Lemmas~\ref{lm:CombinatorialImpliesFibreBijection} and 
Proposition~\ref{prop:FibreOfPrincipalBall} together.


\begin{lm} \label{lm:bijection-of-index-sets}
  Let $U \subset \topo \Delta$ be a principal ball that is open, $V$ in $\conncomp { \inv {\topo \Phi}(U) }$, and $y$ in~$U$. 
  If $\Phi$ is \combinatorial{}, and $\Sigma$ and $\Delta$ have face morphisms that are isometries,
  then $\poly_\Sigma$ induces a bijection 
	\begin{align} \label{eq-bijection-of-index-sets} 
	 V \cap \inv {\topo \Phi} (y)
	 \to \poly_\Sigma (V) \cap \inv \varphi( \poly_\Delta(y) ).
	 \end{align}
\end{lm}

\begin{proof} 
  Since $\Phi$ is \combinatorial{}, by Lemma~\ref{lm:CombinatorialImpliesFibreBijection} the map from Equation~\eqref{eq-bijection-of-index-sets} is injective;
  and for $\alpha$ in $\poly_\Sigma (V) \cap \inv \varphi( \poly_\Delta(y))$, there is $x$ in $\inv {\topo \Phi}(y)$ such that $\poly_\Sigma(x) = \alpha$. 
  Suppose that $x$ is not in $V$. 
  Let $V'$ be the connected component of $x$ in $\inv {\topo \Phi}(U)$. 
  By  Proposition~\ref{prop:FibreOfPrincipalBall},  we have that $\poly_\Sigma(V')$ is a connected component distinct from~$\poly_\Sigma(V)$. 
  But $\alpha$ is both in $\poly_\Sigma(V)$ and $\poly_\Sigma(V')$, a contradiction. 
  Hence, $x$ is in $V$, as desired.
\end{proof}

Having related both index sets, it is straightforward to prove the following:

\begin{prop} 
\label{prop:ReduceToCombi}
    Assume that $\Phi$ is \combinatorial{} and that $\Sigma$ and $\Delta$ have face morphisms that are isometries. 
    Let $m_\varphi : \Sigma \to \ZZ_{\ge 1}$ be an index map.
    We have that: 
    \begin{enumerate} 
      \item If $(\varphi, \, m_\varphi)$ is an indexed branched cover with branch locus $\calB \subset \Delta$,
        then $(\topo \Phi, \, m_\varphi \circ \poly_\Sigma)$ is an indexed branched cover with branch locus $\inv \poly_\Delta \calB$.
      \item If $(\topo \Phi, \, m_\varphi \circ \poly_\Sigma)$ is an indexed branched with branch locus $B \subset \topo \Delta$,
        then $(\varphi, \, m_\varphi)$ is an indexed branched with branch locus $\Delta \setminus \poly_\Delta(\topo \Delta \setminus B)$.
    \end{enumerate}
\end{prop} 

\begin{proof}	
  We first make two observations, and then proceed to prove each item. 

  Observation I: Let $y \in \topo \Delta$ and $\beta \in \Delta$ be such that $\poly_\Delta (y) = \beta$.
  By Lemma~\ref{lm:smallEpsilonImpliesOpen}, if $\varepsilon$ is small enough, then $U_\varepsilon = \upsetballsub \Delta y \varepsilon$ is open, hence $\poly_\Delta (U_\varepsilon) = \upset_\Delta \beta$.
  Let $V$ be a connected component in $\conncomp {\inv {\topo \Phi}(U_\varepsilon) }$.
  By Proposition~\ref{prop:FibreOfPrincipalBall} we have that $\poly_\Sigma(V)$ equals $\upset_\Sigma \alpha$ for some $\alpha$ in $\inv \varphi(\beta)$, and is a connected component of $\inv \varphi(\upset_\Delta \beta)$.
  So Diagram~\ref{dia:super-diagram} restricts to 
  \begin{minipage}{\textwidth}
	\centering
	\vspace{0.7em}
	\[ 
	\begin{tikzcd}
	V \arrow[r,"\poly_\Sigma"] \arrow[d,swap,"\topo \Phi"] &
	\upset_\Sigma \alpha \arrow[d,"\varphi"] \\
	U_\varepsilon \arrow[r,"\poly_\Delta"] & \upset_\Delta \beta.
	\end{tikzcd} 
	\]
	
	\refstepcounter{thm}
	      \label{dia:SuperDiagramRestricted}
	Diagram~\ref{dia:SuperDiagramRestricted}
	\vspace{1em}
	\end{minipage}
  Both vertical maps are continuous. 
  Since $\poly_\Sigma$ and $\poly_\Delta$ are open maps, if one of the vertical maps is an open map, so is the other one. 
  Finally, Lemma~\ref{lm:bijection-of-index-sets} implies that if one of the vertical maps is bijective, so is the other one.   
  Hence, if one of the vertical maps is a homeomorphism, so is the other one.

  Observation II: Let $U \subset \topo \Delta$ be a principal ball that is open, 
  $V$ in $\conncomp { \inv {\topo \Phi}(U) }$, 
  and $y$ a point in $U$. 
  Set $\calU = \poly_\Sigma (U)$. 
  Lemma~\ref{lm:bijection-of-index-sets} then gives
	\begin{align}
    \label{eq:SomeNumerology}
    \deg( {\topo \Phi}, m_\varphi \circ \poly_\Sigma, V)(y) 
    &= \sum_{\substack{ x \in \inv {\topo \Phi}(y) \\ x \in V} } m_\varphi \circ \poly_\Sigma(x) \\
		\nonumber &= \sum_{\substack{ \gamma \in \inv \varphi( \poly_\Delta(y)) \\ \gamma \in \poly_\Sigma(V)} } m_\varphi (\gamma) \\
		\nonumber &= \deg( \varphi, m_\varphi, \poly_\Sigma( V))(\poly_\Delta(y)).
	\end{align}
 
  Proof of Item~(1): Assume that $(\varphi, \, m_\varphi)$ is an indexed branched cover with branch locus $\calB \subset \Delta$, and set $\calW = \Delta \setminus \calB$.
  Since $\poly_\Delta$ is open and continuous, the set $W = \inv \poly_\Delta(\calW)$ is open and dense.
  To see that $\topo \Phi$ is a branched cover with branch locus $ \inv \poly_\Delta \calB = \topo \Sigma \setminus W $,
  let $y$ be in $W$, so $\beta = \poly_\Delta y$ is in $\calW$.
  Consider a principal ball $U_\varepsilon = \upsetballsub \Delta y \varepsilon \subset \topo \Delta$ that is open, so $\poly_\Delta U_\varepsilon = \upset_\Delta \beta$.
  Since $\varphi$ is a branched cover that is unramified over $\calW$ and $\upset_\Delta \beta$ is the smallest open neighbourhood that contains $\beta$, 
  we have that $\inv \varphi(\upset_\Delta \beta)$ comprises several disjoint sets, each mapped to $\upset_\Delta \beta$ homeomorphically by $\dtmor$.
  Since $\upset_\Delta \beta$ is connected, these disjoint sets are the elements of $\conncomp {\inv \varphi(\upset_\Delta \beta)}$. 
  By Observation~I and Proposition~\ref{prop:FibreOfPrincipalBall} this means that each element in $\conncomp {\inv {\topo \Phi}(U_\varepsilon)}$ 
  is mapped homeomorphically to $U_\varepsilon$ by $\topo \Phi$, 
  as desired.

  Since $\topo \Phi$ is a branched cover, to see that $(\topo \Phi, m_\varphi \circ \poly_\Sigma)$ is an indexed branched cover,
  by Lemma~\ref{lm:PathConnectedCountablyBasisSimplifies} and Proposition~\ref{prop:CountableBasisOfPathConnectedSets},
  it is enough to consider
  a principal ball $U = \upsetballsub \Delta y \varepsilon \subset \topo \Delta$ that is open.
  Let $V$ be in $\conncomp {\inv {\topo \Phi}(U)}$.
  By Proposition~\ref{prop:FibreOfPrincipalBall} the open set $\poly_\Sigma(V)$ is a connected component of~$\inv \varphi(\calU)$.
  Since $(\varphi, m_\varphi)$ is an indexed branched cover, 
  $\deg( \varphi, m_\varphi, \poly_\Sigma( V))$ is constant over~$\upset_\Delta \beta$, 
  so Equation~\eqref{eq:SomeNumerology} implies that $\deg( {\topo \Phi}, m_\varphi \circ \poly_\Sigma, V)(y)$ is constant over~$U$.

  Proof of Item~(2): 
  Assume that $(\topo \Phi, \, m_\varphi \circ \poly_\Sigma)$ is an indexed branched cover with branch locus $B \subset \topo \Delta$, and set $W = \topo \Delta \setminus B$.
  Since $\poly_\Delta$ is open and surjective, the set $\calW = \poly_\Delta W$ is open and dense.
  To see that $\varphi$ is a branched cover with branch locus $\Delta \setminus \poly_\Delta(\topo \Delta \setminus B) = \Delta \setminus \calW$, let $\beta$ be in $\calW$.
  Since $\calW$ is open, it is an up-set, we have $\upset_\Delta \beta \subset \calW$, 
  so the set $U = \inv \poly_\Delta( \upset_\Delta \beta)$ is non-empty and contained in $W$.
  Once again we observe that $U$ is connected, so each element in $\conncomp{ \inv {\topo \Phi}(U) }$ gets mapped homeomorphically to $U$ by $\topo \Phi$.
  We conclude again by  Observation~I and Proposition~\ref{prop:FibreOfPrincipalBall}.
   
  Since $\varphi$ is a branched cover, to see that $(\varphi, m_\varphi)$ is an indexed branched cover,  
  by Lemma~\ref{lm:LocalDegreeUnionTwoSets} 
  it is enough to consider
  an upset $\upset_\Delta \beta$. 
  Let $\calV$ be in $\conncomp {\inv \varphi(\upset_\Delta \beta)}$.
  Since $\poly_\Delta$ is surjective, we can choose $y \in \topo \Delta$ such that $\poly_\Delta(y) = \beta$, 
  and by Proposition~\ref{prop:CountableBasisOfPathConnectedSets} we can choose $\varepsilon$ such that $U = \upsetballsub \Delta y \varepsilon$ is open.
  So $\poly_\Delta(U) = \upset_\Delta \beta$.
  By Proposition~\ref{prop:FibreOfPrincipalBall} there is $V$ in $\conncomp{ \inv{\topo \Phi}(U) }$ 
  such that $\poly_\Sigma(V) = \calV$.
  Since $(\topo \Phi, m_\varphi \circ \poly_\Sigma)$ is an indexed branched cover, 
  $\deg( {\topo \Phi}, m_\varphi \circ \poly_\Sigma, V)(y)$ is constant over~$U$. 
  so Equation~\eqref{eq:SomeNumerology} implies that $\deg( \varphi, m_\varphi, \poly_\Sigma( V))$ is constant over~$\upset_\Delta \beta$. \qedhere
\end{proof}

Putting everything together we get Theorem~\ref{thm:SummarySection1}.

\begin{proof}[Proof of Theorem~\ref{thm:SummarySection1}]
The square commutes by Lemma~\ref{lm:superDiagram}.
It is a fibre product by Lemmas~\ref{lm:combinatorialImpliesFibredProduct} and~\ref{lm:phiIsOpen}. 
It remains to apply Proposition~\ref{prop:ReduceToCombi} to get that if one of the vertical maps is an indexed branched cover, so is the other.

Given $\sigma : \Sigma \to \polyintface$, the statement of Theorem~\ref{thm:SummarySection1} is topological, 
so we may replace $\sigma$ with any polyhedral complex whose topological realization is homeomorphic to $\topo \Sigma$, and we may also replace the integral structure.
We apply Example~\ref{ex:HalfSpace} to all the $\alpha \in \Sigma$, so we can assume that the polyhedra of $\Sigma$ are either polytopes or have a unique bounded face of dimension $d-1$, where $d = \dim \alpha$.
After a barycentric subdivision, by Proposition~\ref{prop:BCSIsPolyhedralComplex} and Remark~\ref{re:GeneralizeBarycentricSubdivision}, we get simplices and simplicial cones.
Since Minkowski sums distribute with unions, 
we may assume that every polyhedra $\sigma_\alpha$ of $\Sigma$ is of the form $\conv \calV + \Span_{\RRgo} \calR$ 
where $\calV$ and $\calR$ are subsets of the standard vectors $\aset{e_1, \dots, e_{d}}$ and $d = \dim \sigma_\alpha$.
So we can even assume that $(\sigma_\alpha, N^d) = (\RR_{\ge 0}^{d}, \ZZ^{d})$ and that $f_{\alpha\gamma}$ maps $N^\alpha$ injectively into $N^\gamma$.
This means that we have a canonical euclidean metric $d_\alpha$ on $(\sigma_\alpha, N^\alpha)$ and that any $f_{\alpha\gamma}$ is an isometry. 
Thus, we can now apply Proposition~\ref{prop:ReduceToCombi} to get the theorem. 
\end{proof}

\subsection{Refining morphisms}
    \label{sub:RefiningMorphisms}
Given a morphism $\Phi : \Sigma \to \Delta$ of polyhedral complexes, 
a natural question is whether 
there exists refinements $\gamma : \Sigma' \to \Sigma$ and  $\tau : \Delta' \to \Delta$, and a morphism  $\Phi' : \Sigma' \to \Delta'$ such that 
the topological realizations $\topo \Phi$ and $\topo {\Phi'}$ are isomorphic and $\Phi$ is \combinatorial.
For the case of dimension~1, such result is \emph{Construction~17 (essential model of a tropical morphism)} in \cite{dv20}, exemplified by:

\begin{ex} 
    \label{ex:NotAFiberProduct}
    Consider the morphism $\Phi : [\Gamma : G \to \polyintface] \to [\Delta : T \to \polyintface]$ 
    of one-dimensional polyhedral complexes shown in Figure~\ref{fig:FibreNotEqual}.
    The fibre over $t \in T$ has 3 points, namely $\inv \varphi(t) = \aset{e, B, f}$.
    Observe $\downset_T \dtmor(B) = \aset{u, t, v}$ and $\dtmor(\downset_G B) = \aset{t}$,
    hence $\dtmor$ is not \combinatorial{}.
    Let $z$ be a point in $\Delta_t^\circ$, so $\poly_\mG(z) = t$.
    Set $y = \topo \Phi(B)$.
    The fibre $\inv{ \topo \Phi}(z)$ has 1 point if $z$ is right of $y$, and 2 points otherwise.
    Thus, there is never a bijection between $\inv{ \topo \Phi}(z)$ and $\inv \varphi(t)$.
    Upon a refinement to $\tmor'$, we obtain a \combinatorial{} morphism, 
    and one sees that $\topo {\mG'}$ is the fibre product of $\dtmor'$ and $\poly_T$.
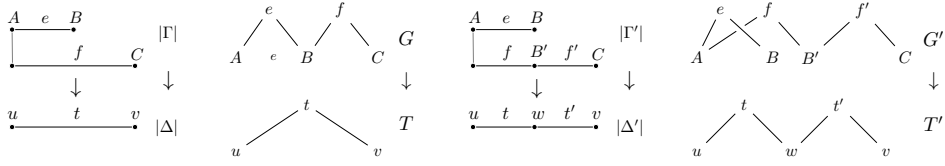
\begin{figure}
\centering

\begin{minipage}{0.20\textwidth}
\scalebox{0.63}{	\begin{overpic}{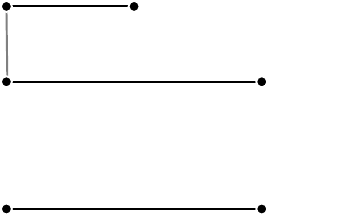}

		\put (-1,66) {\scalebox{1}{$A$}}
		\put (37,66) {\scalebox{1}{$B$}} 
		\put (20,66) {\scalebox{1}{$e$}}
		\put (40,45) {\scalebox{1}{$f$}}
		\put (75,43) {\scalebox{1}{$C$}}
		\put (38,23) {\scalebox{1.3}{$\downarrow$}}
	 
		\put (-1,7) {\scalebox{1.1}{$u$}}
		\put (40,7) {\scalebox{1.1}{$t$}}
		\put (75,7) {\scalebox{1.1}{$v$}}

    \put (95,27) {\scalebox{1.3}{$\downarrow$}}	
		\put (92,57) {\scalebox{1}{$\topo \Gamma$}} 
		\put (90,-1) {\scalebox{1}{$\topo \Delta$}} 
		
\end{overpic}}

\end{minipage}
\begin{minipage}{0.23\textwidth}
	\scalebox{0.63}{
		 \begin{tikzpicture}
		\node (A1) at (0,0) {$A$};
		\node (e) at (0.8,0) { \scalebox{0.8}{$e$} };
		\node (B1) at (1.5,0) {$B$};
		\node (C1) at (3,0) {$C$};
		\node (t1) at (2.2,1) {$f$};
		\node (s2) at (0.7,1) {$e$};
		\draw (B1) -- (t1) -- (C1);  
		\draw (A1) -- (s2) -- (B1);  
		 \node (A) at (0,-2) {$u$};
		 \node (B) at (3,-2) {$v$};
		 \node (s) at (1.5,-1) {$t$};
		 \draw (A) -- (s) -- (B);  
		 \node (Gamma) at (3.6,0.4) {\scalebox{1.2}{$G$}};
		 \node (Arrow) at (3.6,-0.5) {\scalebox{1.2}{$\downarrow$}};
		 \node (t) at (3.6,-1.4) {\scalebox{1.2}{$T$}};
		 \end{tikzpicture}
	}
\end{minipage}
\begin{minipage}{0.20\textwidth}
\scalebox{0.63}{	\begin{overpic}{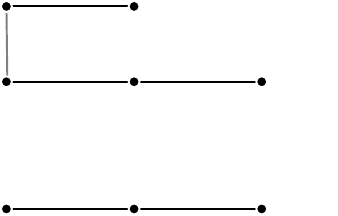} 
		\put (-1,66) {\scalebox{1}{$A$}}
		\put (37,66) {\scalebox{1}{$B$}} 
		\put (37,45) {\scalebox{1}{$B'$}} 
		\put (75,43) {\scalebox{1}{$C$}}
		\put (38,22) {\scalebox{1.3}{$\downarrow$}}
		\put (20,66) {\scalebox{1}{$e$}}
		\put (20,45) {\scalebox{1}{$f$}}
		\put (58,45) {\scalebox{1}{$f'$}}

    \put (-1,7) {\scalebox{1.1}{$u$}}
		\put (20,7) {\scalebox{1.1}{$t$}}
		\put (38,7) {\scalebox{1.1}{$w$}}
		\put (58,7) {\scalebox{1.1}{$t'$}}
		\put (75,7) {\scalebox{1.1}{$v$}}

    \put (95,27) {\scalebox{1.3}{$\downarrow$}}	
    \put (92,57) {\scalebox{1}{$\topo{\Gamma'}$}} 
    \put (90,-1) {\scalebox{1}{$\topo{\Delta'}$}} 

\end{overpic}} 
\end{minipage}
\begin{minipage}{0.23\textwidth}
	\scalebox{0.63}{
		 \begin{tikzpicture}
		\node (A1) at (0,0) {$A$};
		\node (B1) at (1.6,0) {$B$};
		\node (D) at (2.4,0) {$B'$};
		\node (C2) at (4.4,0) {$C$};
		\node (t2) at (3.5,1) {$f'$};
		\node (s1) at (0.5,1) {$e$};
		\node (s2) at (1.5,1) {$f$};
		\draw (A1) -- (s1) -- (B1);  
		\draw (A1) -- (s2) -- (D) -- (t2) -- (C2);  
		 \draw[preaction={draw=white, -,line width=6pt}] (B1) -- (s1);
		 \node (A) at (0,-2) {$u$};
		 \node (B) at (2,-2) {$w$};
		 \node (C) at (4,-2) {$v$};
		 \node (t) at (3,-1) {$t'$};
		 \node (s) at (1,-1) {$t$};
		 \draw (A) -- (s) -- (B) -- (t) -- (C);  
		 \node (Gamma) at (5,0.4) {\scalebox{1.2}{$G'$}};
		 \node (Arrow) at (5,-0.5) {\scalebox{1.2}{$\downarrow$}};
		 \node (t) at (5,-1.4) {\scalebox{1.2}{$T'$}};
		 \end{tikzpicture}
	}
\end{minipage}
	\caption{
  \label{fig:FibreNotEqual}
  On the left, a morphism of 1-dimensional polyhedral complexes together with the corresponding map of posets.
On the right a refinement which makes it a \combinatorial{} map. }
\end{figure}
\end{ex}

Going to higher dimensions,
the pushforward of the polyhedral structure of the source $\Sigma$ no longer induces a polyhedral complex on the target $\Delta$, 
see Example~\ref{ex:ComplicatedRefinement}.
Further refinement is needed via stellar subdivisions, which we exemplify now in the embedded case.
The generalization to polyhedral spaces is straightforward;  see \cite[Section 9.9]{var22} for a detailed description of the face poset and the functor to $\coneintface$.

\begin{figure} 
  \centering

    	\noindent \begin{minipage}{0.2\textwidth}
		  \centering
		      \begin{overpic}[scale=0.5]{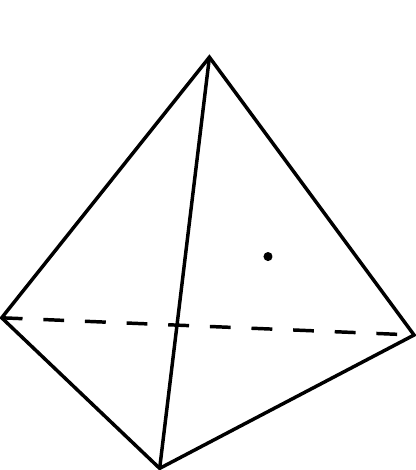}
		      \end{overpic}  
	    \end{minipage}
	    \begin{minipage}{0.2\textwidth}
		  \centering
		      \begin{overpic}[scale=0.5]{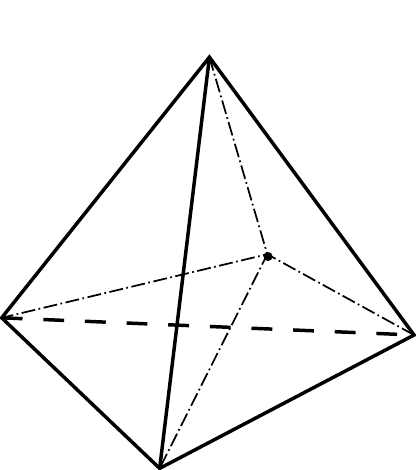}
		      \end{overpic}
	    \end{minipage}
      
      \caption{
        \label{fig:StellarSubdivision} 
        On the left, the intersection of the positive orthant of $\RR^4$ with the hyperplane $\sum x_i = 1$,
        and $p = (1/3,1/3,1/3,0)$. 
        On the right, the stellar subdivision by $p$,
        which adds 1 cone of dimension~1, 4 cones of dimension~2, 6 cones of dimension~3 and 3 cones of dimension~4.
      }
  \end{figure}

\begin{ex}
      \label{ex:StellarSubdivisionOfACone}
  Let $\sigma$ be a rational polyhedral cone in $(V,N)$ and $x$ a point in $\sigma$.
  Let $\PolyStar \sigma x$ be the set of faces of $\sigma$ that contain $x$.
  The stellar subdivision 
  \[\SD \sigma x = \aset{ \Span_{\RRgo}(\tau, x) \suchthat \tau \text{ face of } \sigma, \,  \tau \not \in \PolyStar \sigma x   }\]
  is an embedded polyhedral complex
  and the morphism $\SD \sigma x \to \sigma$ given by inclusion is a refinement.
  See Figure~\ref{fig:StellarSubdivision} for a subdivision of a 4-dimensional simplicial cone.  
  \end{ex}

\begin{ex} 
      \label{ex:StellarSubdivisionOfAComplex}
  Consider an embedded cone complex $\Sigma$ in $(V, N)$, and a point $x \in \topo \Sigma$.
  We get a refinement $\SD \Sigma x \to \Sigma$ by performing a stellar subdivision on each cone containing~$x$.
  Concretely, 
  if we let $\delta \in \Sigma$ be the unique cone with $x \in \delta^\circ$, 
  the cones deleted are those in 
  $\PolyStar \Sigma \delta = \aset{ \sigma \in \Sigma \suchthat \delta \subset \sigma  }$.
  and a cone $\Span_{\RRgo}(x, \tau)$ is added for each $\tau$ not-containing $x$ that is a face of a cone $\sigma$ in $\PolyStar \Sigma \delta$. 
  That is,
  \[ \SD \Sigma x = 
   ( \Sigma \setminus  \PolyStar \Sigma \delta ) \cup
  \aset{\Span_{\RRgo}(\tau, x) \suchthat \tau \not \in \PolyStar \Sigma \delta \text{ and } \tau \subset \sigma \text{ for some } \sigma \in \PolyStar \Sigma \delta  }.\]
  See \cite[Definition~2.1]{ewa12}  for more information, 
  including the generalization to the cell complex case.
  See \cite[Definition 2.22]{koz07} for the abstract simplicial complex case.
\end{ex}

On the other hand, the pullback refinement is accomplished via half space subdivisions. 

\begin{ex}
      \label{ex:HalfSpaceSubdivisionOfAPolyhedron}
  Let $\sigma$ be a rational polyhedron in $(V,N)$, again we assume for a moment that $\Span \sigma$ is not necessarily the whole $V$.
  Given a rational functional $u \in N^*$ we get a subdivision of $\sigma$ by considering the two polyhedra $\sigma^+ = \sigma \cap H^+(u,c)$ and $\sigma^- = \sigma \cap H^-(u,c)$,
  where $H^+$ and $H^-$ are the half-spaces introduced at the beginning of Subsection~\ref{sub:PolyhedralSpaces}.
\end{ex}

\begin{cons} 
    \label{cons:PullBackConstruction}
Let $F : \topo \Sigma \to \topo \Delta$ be a continuous map that is locally affine and injective on each $\sigma \in \Sigma$.
Given $\sigma \in \Sigma$, consider $\delta_1, \dots, \delta_\mu$ the polyhedra in $\poly_\Delta(F(\sigma))$ whose dimension is $\dim \sigma -1$.
Since $F$ is injective on $\sigma$,
if there are no such polyhedra, $F(\sigma) \subset \delta^\circ$ for some $\delta \in \Delta$ and there is no subdivision to make.
Otherwise, do a half-space subdivision with the $\RR$-span in $N^\sigma_\RR$ of each $\inv F(\delta_i)$, as in Example~\ref{ex:HalfSpaceSubdivisionOfAPolyhedron}.
Repeating this process for all $\sigma \in \Sigma$ yields the \emph{pullback refinement} $F^* : \Sigma^* \to  \Sigma$.
\end{cons}

\begin{lm} 
    \label{lm:PullbackIsNice}
If $F : \topo \Sigma \to \topo \Delta$ is a continuous map between  polyhedral complexes such that for each $\sigma \in \Sigma$ the resctriction $F|_{\topo \sigma}$ is an affine injective map, then the pullback refinement $\Sigma^*$ is a polyhedral complex.
\end{lm}

We refine both $\Sigma$ and $\Delta$ by alternating stellar subdivisions 
and pullback refinements. 
For ease of exposition, and also eyeeing on future applications to tropical moduli spaces, we assume that $\Sigma$ and $\Delta$ are spaces of cones. 
The generalization is straightforward, following the lines laid down in Remark~\ref{re:GeneralizeBarycentricSubdivision}.
The following construction can be best understood while reading along Example~\ref{ex:ComplicatedRefinement}.

\begin{cons} 
    \label{cons:WonderfulConstruction}
    Let $F : \topo {\Sigma_0} \to \topo {\Delta_0}$ be a continuous map that is locally affine and injective on each $\sigma \in \Sigma$.
Let $\tau_1, \dots, \tau_\mu$ be the set of rays of $\Sigma_0$ ordered in such a way that 
\begin{align} \label{eq:DimCondition} \dim \poly_\Delta F(\tau_i) \le  \dim \poly_\Delta F(\tau_j) \end{align}
whenever $i \le j$. 
Iterate over the $\tau_i$ doing the following:
\begin{enumerate} 
  \item Do a stellar subdivision on $\Delta$ centered at $F(\tau_i)$ to get $\Delta_{i+1}$
  \item Pullback the new polyhedral structure to obtain $\Sigma_{i+1}$.
  \item If there are new rays in $\Sigma_{i+1}$, relabel the ray sequence such that the condition from Equation~\eqref{eq:DimCondition} still holds, and also the newer rays of $\Sigma_{i+1}$ appear after the original rays of $\Sigma_i$ of the same dimension.
  \item Repeat until there are no more rays to subdivide from, that is until $F(\tau_i)$ is a ray of $\Delta_k$ for all $i$.
\end{enumerate}
\end{cons}

The proof that the construction above terminates is analogous to the usual proof that the barycentric subdivision can be obtained by doing stellar subdivisions on the barycentres. 
The key insight is the way that the rays have been ordered. 
Unlike the case of the BCS, here the final complex does depend on the chosen order for the rays. 
 
\begin{ex} 
    \label{ex:ComplicatedRefinement}
   The following sequence of figures illustrate Construction~\ref{cons:WonderfulConstruction}.

\vspace{1em}

	\noindent
  \begin{minipage}[t]{\textwidth}
	\begin{minipage}[t]{.3\textwidth}

		\begin{overpic}[scale=1.2]{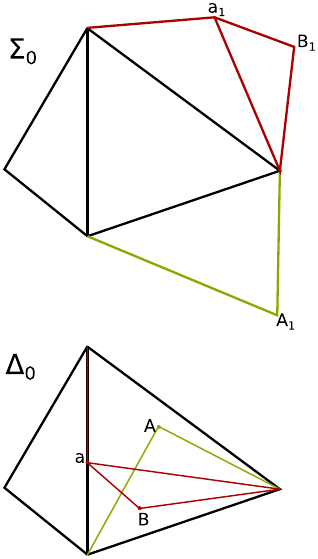} 
	\end{overpic}
  
\center
		\vspace{0.5em}

    (a)	
		
	\end{minipage}\hspace{1em}
	\noindent \begin{minipage}[t]{.3\textwidth} 
		\begin{overpic}[scale=1.2]{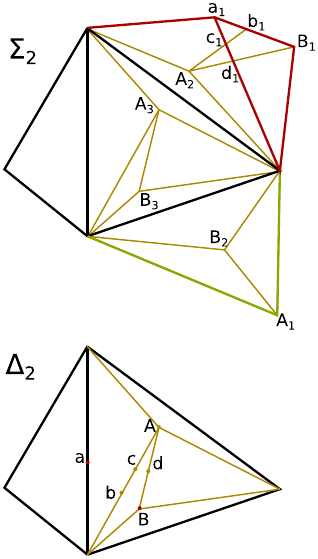} 
	\end{overpic} 

	\center	
		\vspace{0.5em}

  (b)

	\end{minipage}\hspace{1em}
	\begin{minipage}[t]{.3\textwidth}
		\begin{overpic}[scale=1.2]{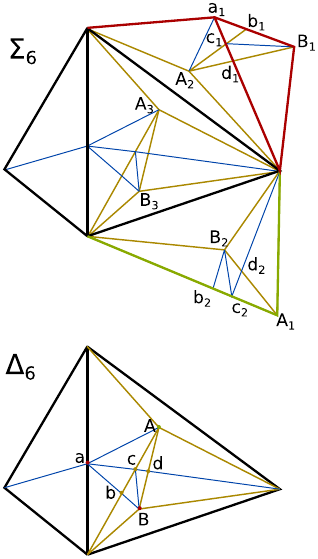} 
	\end{overpic} 
	\center	

		\vspace{0.5em}
  (c)
		
	\end{minipage}\hspace{1em}
\vspace{1em}
\refstepcounter{figure}
\label{fig:CombinatorialSubdivision}
{\sc Figure~\ref{fig:CombinatorialSubdivision}}. On the left (a), the initial map $\Sigma_0 \to \Delta_0$. 
  Note that the projection does not define a polyhedral complex on $\Delta_0$.
  In the centre (b), stellar subdivision on $\Delta$, first $A$ and then $B$, and pullbacks to~$\Sigma$ that produce $b$, $c$, and~$d$ as images of rays.
  On the right (c), stellar subdivision of $a$, $b$, $c$, $d$ on $\Delta_2$ and pullbacks to~$\Sigma_2$.

  \end{minipage}

\end{ex}

\begin{re}[induced affine map on polyhedral spaces] 
    \label{re:Barycentric}
    Let $\sigma : \Sigma \to \polyintface$ and $\delta : \Delta \to \polyintface$ be polyhedral complexs, thus potentially with non-trivial automorphisms and stack behaviour. 
    Let $\Phi : [\Sigma \to \polyintface ] \to [\Delta \to \polyintface] $ be a morphism.
    This morphism induces a continuous map that is locally affine and locally injective on the barycentric subdivisions, namely an $F : \BCS \Sigma \to \BCS \Delta$ to which Construction~\ref{cons:WonderfulConstruction} can be applied.  
\end{re}




 
           \section{The balancing condition and connectivity} 
                \label{sec:TheBalancingConditionAndConnectivity}


\mysubsection{The balancing condition}
    \label{TheBalancingCondition}
    The usual definition of a tropical morphism of metric graphs, i.e. of 1-dimensional polyhedral complexes,
    starts with a continuous map that is piecewise-linear with integral slopes, i.e. a morphism of polyhedral complexes.
    This map must satisfy two \emph{realizability conditions} that brings it closer to the theory of Riemann Surfaces:
    the balancing condition and the Riemann-Hurwitz condition. 
    Now we focus on the former condition, to show that it is equivalent to the map being an indexed branched cover, and this proof works for higher dimensions.

 \begin{de}
	Let $\varphi : \Sigma \to \Delta$ be a morphism of finite posets, 
  $\calV \subset \Sigma$ an up-set, 
  and $m_\calV : \calV \to \ZZ_{\ge 1}$ a map. 
  We say that $(\varphi, m_\calV)$ is \emph{balanced},
  if for any $\alpha$ in $\calV$ and any choice of $\beta$ in $\Delta$ such that $\varphi(\alpha) \lessdot \beta$
  the following equation, called \emph{the balancing condition}, holds: 
  \begin{align} \label{eq:balancing-condition}
  m_\calV(\alpha) = \sum_{ \substack{ \gamma \in \inv \dtmor(\beta) \\ \alpha \lessdot \gamma} } m_\calV(\gamma).
 \end{align}
\end{de}

\begin{ex}[discrete tropical morphisms] \label{re:discrete-tropical-morphism}
    	Recall that using the index map from Example~\ref{ex:TropicalMorphism}, 
      the map $\mG \to \mT$ of posets from Figure~\ref{fig:covers}~(b) becomes an indexed branched cover. 
      One can verify the balancing condition.
 
      By Theorem~\ref{thm:SummarySection1}, the topological realization $\topo \Phi : \topo \mG \to \topo \mT$ 
      is an indexed branched cover of metric graphs. 
      Moreover, since the target of $\topo \Phi$ is simply connected, 
      $\mG$ can be regarded as several copies of $\topo \Delta$ glued at certain regions. 
      This is pictured in Figure~\ref{fig:covers}~(c), where the dashed lines represent identification of points, as done in \cite{dv20}.
\end{ex}

If $\varphi$ is clear from context,
we say that $m_\calV$ is a \emph{balanced map}. 
The convention for empty sums in Equation~\eqref{eq:balancing-condition} is that they evaluate to zero.
This implies that if $m_\calV$ is a balanced map, and $\alpha$ is in $\max \calV$, then we must have that $\varphi(\alpha)$ is in $\max \Delta$,
for otherwise  we could take a $\beta$ covering $\varphi(\alpha)$ and apply the balancing condition to get an element in $\inv \varphi(\beta)$ that covers $\alpha$.
Moreover, if $\varphi$ is \combinatorial{} we get:

\begin{lm} 
        \label{lm:BalancedImpliesBranchedCover}
      Let $\varphi : \Sigma \to \Delta$ be a \combinatorial{} morphism of posets,  
      $\calV \subset \Sigma$ an up-set, 
      and $m_\calV : \calV \to \ZZ_{\ge 1} $ a balanced map.
      The corestriction $\varphi \corestrictedto {\varphi(\calV)} : \calV \to \varphi(\calV)$ 
      is a branched cover with branch locus $\varphi(\calV) \setminus \max \Delta$.
\end{lm}

\begin{proof} 
      We check the two conditions from Example~\ref{ex:BranchedCoverPosets}.
      To check the first condition, let $\beta$ be in $\max \varphi(\calV)$.
      Suppose there is $\alpha \in \inv{\varphi}(\beta)$ such that $\alpha$ is not maximal, 
      i.e.~there is $\gamma$ in $\Sigma$ such that $\alpha < \gamma$.
      Since $\varphi$ is \combinatorial{}, it maps $\downset_\Sigma \gamma$ isomorphically to $\downset_\Delta \varphi(\gamma)$,
      thus $\varphi(\gamma)$ is strictly greater than $\beta$, contradicting that $\beta$ is maximal.
      Hence, $\inv{\varphi}(\beta) \subset \max \Sigma$. 
      The second condition is true by definition. 
\end{proof}

Thus, given two elements $\mu \lessdot \nu$, we investigate the fibres above them.

\begin{lm}
  \label{lm:TreeStructure} 
	Let $\varphi : \Sigma \to \Delta$ be a \combinatorial{} morphism of posets.
  For any $\alpha$ in~$\Sigma$, and $\mu, \nu$ in $\upset_\Delta \varphi(\alpha)$ such that $\mu \lessdot \nu$, we have the following partition: 
\begin{align}
  \label{eq:PartitionMu}
      \inv{\varphi}(\nu) \cap  \upset_\Sigma \alpha = \bigsqcup_{ \substack{
                                                       \gamma \in \inv{\varphi}(\mu) \\ 
                                                       \gamma \in \upset_\Sigma \alpha}}
                                                                  \inv{\varphi}(\nu) \cap \aset{\eta \in \Sigma \suchthat \gamma \lessdot \eta }.   
\end{align}
\end{lm}
\begin{proof}
  Consider the family $\aset{\upset_\Sigma \gamma }$ indexed by $\gamma \in \inv{\varphi}(\mu) \cap \upset_{\Sigma} \alpha$. 
  Each member of $\aset{\upset_\Sigma \gamma}$ is a subset of~$\upset_\Sigma \alpha$, and by Lemma~\ref{lm:connected-components-combinatorial-morphism} they are pairwise disjoint. 
  Since $\aset{\eta \in \Sigma \suchthat \gamma \lessdot \eta } \subseteq \upset_\Sigma \gamma$,  
  we are done if we show that the right hand side of Equation~\eqref{eq:PartitionMu} contains the left hand side. 
  Namely, for $\eta$ in $\inv{\varphi}(\nu) \cap \upset_\Sigma \alpha$, show there is $\gamma \in \inv{\varphi}(\mu) \cap \upset_\Sigma \alpha$ such that $\gamma \lessdot \eta$. 
  Since $\varphi$ is \combinatorial{}, it maps $\downset_\Sigma \eta$ isomorphically to $\downset_\Delta \varphi(\eta)$, so there is a unique $\gamma$ such that $\gamma$ is in $\downset_\Sigma \eta$ and $\varphi(\gamma) = \mu$. Since  $\varphi(\eta) = \nu$ and $\mu \lessdot \nu$, we have that $\gamma \lessdot \eta$. Moreover, $\upset_\Sigma \gamma$ is a connected set that intersects at $\eta$ the set $\upset_\Sigma \alpha$. The latter is a connected component of the space $\inv{\varphi}(\upset_\Delta \varphi(\alpha) )$, thus $\upset_\Sigma \gamma \subset \upset_\Sigma \alpha$, so $\gamma \in \upset_\Sigma \alpha$ as desired.   
\end{proof}  

\begin{re} 
    \label{re:TreeStructure}
    Recall that for a given poset $\Sigma$, the \emph{covering graph} $\tG_\Sigma$  has as vertices the elements of $\Sigma$ and as edges the covering relations $\alpha \lessdot \gamma$. 
    Intuitively, what Lemma~\ref{lm:TreeStructure} says is that for any upwards path $P : \varphi(\alpha) \lessdot \beta_1 \lessdot \dots \lessdot \beta_k$, 
    the preimage $\inv \varphi(P)$ has a tree structure.
    In other words, $\tG_{\inv \dtmor(P)}$ is a forest, namely a disjoint union of trees.
\end{re}

The above-proven tree structure implies some crucial formulas.

\begin{lm} 
  \label{lm:BalancingFormulas}
   Let $\varphi : \Sigma \to \Delta$ be a \combinatorial{} morphism of posets,
   $\calV$ an up-set,
   and $m_\calV : \calV \to \ZZ_{\ge 1} $ a balanced map. 
   For any $\alpha$ in~$\Sigma$, 
   and $\mu, \nu$ in $\upset_\Delta \varphi(\alpha)$ 
   such that $\mu \lessdot \nu$
   and $\inv{\varphi}(\mu) \cap \upset_\Sigma \alpha$ is contained in $\calV$, 
   we have:
 \begin{align}
   \label{eq:MovingBalancingAround}
      \sum_{\substack{\gamma \in \inv{\varphi}(\mu) \\ \gamma \in \upset_\Sigma \alpha}} m_\calV(\gamma)
       =  \sum_{\substack{\gamma \in \inv{\varphi}(\mu) \\ \gamma \in \upset_\Sigma \alpha}} \,\, \sum_{\substack{\eta \in \inv{\varphi}(\nu)  \\ \gamma \lessdot \eta }}  m_\calV(\eta)
      = \sum_{\substack{\eta \in \inv{\varphi}(\nu) \\ \eta \in \upset_\Sigma \alpha}} m_\calV(\eta) 
\end{align}
    In particular, if $\beta \in \Delta$ is such that $\varphi(\alpha) \lessdot \beta$, we have:
    \begin{align} \label{eq:BranchedCoverVSBalancingCondition}
      \sum_{\substack{\eta \in \inv{\varphi}(\beta)\\ \eta \in \upset_\Sigma \alpha }}m_\calV(\eta )= \sum_{ \substack{ \eta \in \inv \dtmor(\beta) \\ \alpha \lessdot \eta} } m_\calV(\eta).
    \end{align}
\end{lm}

\begin{proof} 
   The first equality in Equation~\eqref{eq:MovingBalancingAround} is implied by the balancing condition. 
   By Lemma~\ref{lm:TreeStructure} the index sets of the sums in the middle constitute a partition of the index sets of the sums on the right hand side, so the second equality folows.
   We obtain Equation~\eqref{eq:BranchedCoverVSBalancingCondition} from Equation~\eqref{eq:MovingBalancingAround} by setting $\mu = \varphi(\alpha)$, $\nu = \beta$ 
   and noting that the index set for the sum on the left becomes  $\inv{\varphi}(\alpha) \cap \upset_\Sigma(\alpha) = \aset{\alpha}$ a singleton, 
   because $\varphi$ is \combinatorial{}.
\end{proof}

   \begin{re}
      Note that by Equation~\eqref{eq:degree-fmv},
      the left hand side in Equation~\eqref{eq:MovingBalancingAround} equals $\deg(\varphi, m_\varphi, \upset_\Sigma \alpha)(\mu)$ 
      and the right hand side equals $\deg(\varphi, m_\varphi, \upset_\Sigma \alpha)(\nu)$. 
   \end{re}

\begin{re} 
    \label{re:alphaNotInImageOfCalV}
    Crucially, note that in Lemma~\ref{lm:BalancingFormulas} we do not need $\alpha$ to be in the image $\varphi(\calV)$.
\end{re}

   Applying the previous formulas, we arrive to a characterization of indexed branched covers of posets. 
   
	\begin{prop}
     \label{prop:IndexedBranchedCovervsBalanced} 
		Let $\varphi : \Sigma \to \Delta$ be a \combinatorial{} morphism of posets,
    and $m_\varphi : \Sigma \to \ZZ_{\ge 1}$ an index map. We have that
    \begin{itemize} 
      \item If $(\varphi, m_\varphi)$ is an indexed branched cover, then $(\varphi, m_\varphi)$ is balanced.
      \item If $m_\varphi$ is balanced, then $\varphi \corestrictedto{\varphi(\Sigma)} :  \Sigma \to \varphi(\Sigma)$ is an indexed branched cover.
    \end{itemize}
	\end{prop} 

	\begin{proof}
    Let $\alpha$ be in $\Sigma$, $\beta$ in  $\Delta$ with $\varphi(\alpha) \lessdot \beta$, and $\calU = \upset _{\Delta} \varphi (\alpha)$.
    By Lemma~\ref{lm:connected-components-combinatorial-morphism} the connected component $\calV$ of
    $\inv \varphi(\calU)$ that contains $\alpha$ is the upper-set $\upset_{\Sigma} \alpha$.
    So if $(\varphi, m_\varphi)$ is an indexed branched cover, by Equation~\eqref{eq:BranchedCoverVSBalancingCondition} we have
    \begin{align*} 
      m_\varphi(\alpha) = \sum_{\substack{\eta \in \inv{\varphi}(\beta)\\ \eta \in \upset_\Sigma \alpha }}m_\varphi(\eta )= \sum_{ \substack{ \eta \in \inv \dtmor(\beta) \\ \alpha \lessdot \eta} } m_\varphi(\eta).
    \end{align*}

    Now assume that $m_{\varphi}$ is balanced. 
    By Lemma~\ref{lm:BalancedImpliesBranchedCover} the corestriction $\tilde \varphi = \varphi \corestrictedto{\varphi(\Sigma)}$ 
    is a branched cover with branch locus equal to $\varphi(\Sigma) \setminus \max \Delta$.
    By Lemmas~\ref{lm:LocalDegreeUnionTwoSets} and~\ref{lm:connected-components-combinatorial-morphism}, 
    we show that $\deg(\tilde \varphi, m_\varphi, \upset_\Sigma \alpha)$ is constant over $\upset_\Delta \tilde \varphi(\alpha)$. 
    This is clear from Equation~\eqref{eq:MovingBalancingAround} since the left hand side is $\deg(\tilde \varphi, m_\varphi, \upset_\Sigma \alpha)(\mu)$  and the right hand side is $\deg(\tilde \varphi, m_\varphi, \upset_\Sigma \alpha)(\nu)$. 
    So we are done by applying over a sequence that goes from $\tilde \varphi(\alpha)$ up to any desired $\nu$ in $\upset_\Delta \tilde \varphi(\alpha)$, such that each element covers the previous one. 
  \end{proof}

\begin{re} 
    \label{re:CounterexamplesToIBCiffBalanced}
       In Proposition~\ref{prop:IndexedBranchedCovervsBalanced}, the condition on $\varphi$ being \combinatorial{} is crucial.
       Figure~\ref{fig:CounterexamplesToIBCiffBalanced} shows two possible failures when this condition is absent.
 \begin{figure} \center
 
    \begin{tikzpicture}
       \node (A1) at (0,0) {\scalebox{0.85}{$(A_1,1)$}};
       \node (A2) at (2,0) {\scalebox{0.85}{$(A_2,1)$}};
       \node (B1) at (1,1.3) {\scalebox{0.85}{$(B_1,1)$}};
       \node (A) at (3,0) {\scalebox{0.85}{$A$}};
       \node (B) at (3,1.3) {\scalebox{0.85}{$B$}};
    \draw (A1) -- (B1) -- (A2);
    \draw (A) -- (B);  
    \node (varphi1) at (2.3, -0.6) {$\varphi_1 : \Sigma_1 \to \Delta_1$};
 	\end{tikzpicture}
\hspace{4em}
 \begin{tikzpicture}
       \node (A1) at (0,0) {\scalebox{0.85}{$(A_1,2)$}};
       \node (A2) at (2,0) {\scalebox{0.85}{$(A_2,2)$}};
       \node (B1) at (-0.5,1.3) {\scalebox{0.85}{$(B_1,1)$}};
       \node (B2) at (1,1.3) {\scalebox{0.85}{$(B_2,2)$}};
       \node (B3) at (2.5,1.3) {\scalebox{0.85}{$(B_3,1)$}};
       \node (A) at (4,0) {\scalebox{0.85}{$A$}};
       \node (B) at (4,1.3) {\scalebox{0.85}{$B$}};
       \draw (B1) -- (A1) -- (B2) -- (A2) -- (B3);
    \draw (A) -- (B);  
    \node (varphi1) at (3, -0.6) {$\varphi_2 : \Sigma_2 \to \Delta_2$};
 	\end{tikzpicture}

\caption{\label{fig:CounterexamplesToIBCiffBalanced} 
Two morphisms of posets $\varphi_1, \varphi_2$ and maps $m_1 : \Sigma_1 \to \ZZ_{\ge 1}$, $\, m_2: \Sigma_2 \to \ZZ_{\ge 1}$ given by the second numbers in the pairs in the diagrams;
$(\varphi_1, m_1)$ is balanced
but not an indexed branched cover,
and $(\varphi_2, m_2)$ is an indexed branched cover
but not balanced. 
}
\end{figure}
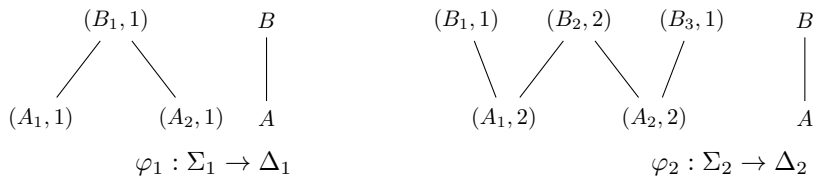      
On the left we have that $(\varphi_1, m_1)$ is balanced but is not an indexed branched cover, since the count of points in $\inv{\varphi_1}(A)$ is 2 and in $\inv{\varphi_1}(B)$ is 1.
On the right, the preimage count is right, but neither  $A_1$ nor $A_2$ satisfy the balancing condition, e.g.~the value of $m_2(A_1)$ would need to be 3.
\end{re}

Note that if $\varphi$ equals its corestriction $\varphi \corestrictedto {\varphi(\Sigma)}$, 
then Proposition~\ref{prop:IndexedBranchedCovervsBalanced} implies Theorem~\ref{mytheorem:IndexedBranchedCovervsBalanced}.
The missing ingredient is Proposition~\ref{myprop:IndexedBranchedCoversAreOpenMaps}, which we prove in the next section.

\mysubsection{A necessary condition for Question~\ref{que:Criteria}}  
    \label{sub:ANecessaryConditionForQuestion}
We prove Proposition~\ref{myprop:IndexedBranchedCoversAreOpenMaps} by studying how to lift paths from $\Delta$ to $\Sigma$ under a morphism $\varphi : \Sigma \to \Delta$. 

\begin{de} 
  Given a poset $\Sigma$, its \emph{comparability graph} $G_\Sigma$
  has as vertices the elements of $\Sigma$, 
  and two vertices are joined by an edge if and only if they are comparable.
\end{de}

Note that $G_\Sigma$ is a simple graph, 
so a path in $G_\Sigma$ is a sequence $\langle \gamma_0, \gamma_1, \dots, \gamma_q \rangle$
of elements in $\Sigma$ such that consecutive elements are comparable,
i.e.~$\gamma_{i-1} \le \gamma_i$ or $\gamma_{i-1} \ge \gamma_i$ for $1 \le i \le q$.
A graph is connected if its vertices are pairwise connected by a path.

\begin{lm} 
  \label{lm:ConnectedImpliesSequence}
Let $\Sigma$ be a finite poset with the poset topology.
Then $\Sigma$ is connected as a topological space
if and only if the comparability graph $G_\Sigma$ is connected as a graph. $\qed$
\end{lm}

\begin{re}
  \label{re:ConnectednessOfComparability}
  If $\calV \subset \Sigma$ is connected, we can consider the induced poset on $\calV$ and apply Lemma~\ref{lm:ConnectedImpliesSequence} 
  to conclude the existence of a sequence inside of $\calV$ with consecutive elements that are comparable, between any two elements of $\calV$.
\end{re}




We now use the balancing condition to lift paths.
Given a path $P = \langle \beta_0, \beta_1, \dots, \beta_q \rangle$ in~$\Delta$,
a \emph{lift} is a path $\tilde P = \langle \gamma_0, \gamma_1, \dots, \gamma_q \rangle$ in $\Sigma$ 
such that $\varphi(\gamma_i) = \beta_i$ for $i \in \aset{0,\dots,q}$.

\begin{lm} 
    \label{lm:LiftingUpwardsPaths}  
      Let $\varphi : \Sigma \to \Delta$ be a morphism of posets,
      $\calV \subset \Sigma$ an up-set,
      and $\alpha$ in $\calV$.
      Assume there is a balanced map $m_\calV$.
      Any upwards path $\beta_0 < \beta_1 < \dots < \beta_k$ in $\upset_\Delta \varphi(\calV)$
      with $\beta_0 = \varphi(\alpha)$
      lifts to a path $\tilde P = \langle \gamma_0, \gamma_1, \dots, \gamma_q \rangle$ in $\calV$ with $\gamma_0 = \alpha$.
\end{lm}

\begin{proof} 
      Proceeding by induction, suppose that $\beta_i$ lifts to $\gamma_i$ in $\calV$.
      We search a lift $\gamma_{i+1} \in \calV$ of $\beta_{i+1}$, 
      that is comparable to $\gamma_i$.
      Since $\beta_i \le \beta_{i+1}$ and $\Delta$ is finite, 
      we can choose a sequence $\mu_0 \lessdot \mu_1 \lessdot \dots \lessdot \mu_k$ in $\Delta$
      with $\mu_0 = \beta_i$ and $\mu_k = \beta_{i+1}$.
      Note that $\gamma_i$ is a lift of $\mu_0 = \beta_i$,
      and that $\gamma_i$ is in $\calV$ so $m_\calV(\gamma_i) \ge 1$,
      Hence, successive applications of the balancing condition 
      give lifts $\eta_0 \lessdot \eta_1 \lessdot \dots \lessdot \eta_k = \beta_{i+1}$ of the $\mu_j$, and we set $\gamma_{i+1} = \mu_k$.
      Since $\gamma_i \le \gamma_{i+1}$ and $\calV$ is an up-set, we have $\gamma_{i+1} \in \calV$.
  \end{proof}

  As a corollary of Lemma~\ref{lm:LiftingUpwardsPaths} we get that $\varphi(\calV) = \upset_\Delta \varphi(\calV)$, i.e.~$\varphi(\calV)$ is open.
  Indeed, if $\beta$ is in $\upset_\Delta \varphi(\calV)$, then there is $\mu = \varphi(\alpha)$ in $\varphi(\calV)$ such that $\mu \le \beta$.
  If $\mu = \beta$ we are done, otherwise we can lift the upwards path $\mu < \beta$ to a path $\alpha < \gamma$ with $\gamma$ in $\calV$ and $\beta = \varphi(\gamma)$, 
  showing that $\beta$ is in $\varphi(\calV)$.
  In particular, this proves Proposition~\ref{myprop:IndexedBranchedCoversAreOpenMaps} 

\begin{proof}[Proof of Proposition~\ref{myprop:IndexedBranchedCoversAreOpenMaps}]
     Let $\calV$ be an up-set.  
     The restriction $m_\calV = m_\Sigma \restrictedto \calV$ is a balanced map, 
     because for any element $\alpha \in \calV$ the elements that cover $\alpha$ are also in $\calV$.
     By Lemma~\ref{lm:LiftingUpwardsPaths} and the discussion after it, $\varphi(\calV)$ is an upset.
     Hence $\varphi$ is an open map.
\end{proof}

\begin{re} 
  Combining Lemma~\ref{myprop:IndexedBranchedCoversAreOpenMaps} and Proposition~\ref{prop:IndexedBranchedCovervsBalanced}
  we get that if $\varphi : \Sigma \to \Delta$ is a \combinatorial{} morphism of posets, 
  then  $\varphi$ being an open map is a necessary condition 
  for the existence of an index map $m_\varphi$ such that $(\varphi, m_\varphi)$ 
  is an indexed branched cover.
  This criterion rules out the existence 
  of an index map making the morphism from  Example~\ref{ex:TooManyComponents} an indexed branched cover. 
\end{re}

  Since \combinatorial{} morphisms of posets are closed maps, we can now finish the proof of Theorem~\ref{mytheorem:IndexedBranchedCovervsBalanced}

  \begin{proof}[Proof of Theorem~\ref{mytheorem:IndexedBranchedCovervsBalanced}]
     As remarked earlier, we are done by Proposition~\ref{prop:IndexedBranchedCovervsBalanced} if we show that $\dtmor$ is surjective.
     By Proposition~\ref{myprop:IndexedBranchedCoversAreOpenMaps} we have that $\dtmor$ is an open map.
     Since $\dtmor$ is \combinatorial{}, closed sets $\downset_\Sigma \alpha$ are mapped to closed sets, 
     and since every closed set $\calC$ in $\Sigma$ is a finite union of $\bigcup_{\alpha \in \calC} \downset_\Sigma \alpha$, we have that $\dtmor$ is a closed map as well.
     Thus, $\dtmor(\Sigma)$ is both an open and a closed set in $\Delta$.
     Since $\Delta$ is connected, the only such set is $\Delta$ itself. 
  \end{proof}

\mysubsection{Extending balanced maps} 
     \label{sub:ExtendingBalancedMaps}
  As the index map in Proposition~\ref{prop:IndexedBranchedCovervsBalanced} is defined over the whole domain,
  this raises the question of what can be said when dealing with a balanced map $m_\calV$ defined over a proper subset $\calV \subsetneq \Sigma$. 

\begin{figure} \center
  \scalebox{0.8}{
    \begin{tikzpicture}
       \node (A1) at (-0.7,0) {\scalebox{0.85}{$(A_1,3)$}};
       \node (B1) at (2,0) {\scalebox{0.85}{$(B_1,2)$}};
       \node (B2) at (3,0) {\scalebox{0.85}{$(B_2,1)$}};
       \node (C1) at (5,0) {\scalebox{0.85}{$(C_1,1)$}};
       \node (C2) at (6,0) {\scalebox{0.85}{$(C_2,2)$}};
       \node (beta1) at (2,1.8) {\scalebox{0.85}{$(\beta_1,2)$}};
       \node (beta2) at (3,1.8) {\scalebox{0.85}{$(\beta_2,1)$}};
       \node (gamma1) at (5,1.8) {\scalebox{0.85}{$(\gamma_1,1)$}};
       \node (gamma2) at (6,1.5) {\scalebox{0.85}{$(\gamma_2,2)$}};
    \node (O1) at (1,-1.5) {$O_1$};
    \node (tildeO1) at (4,-1.5) {$\tilde O_1$};
    \node (tildeO2) at (5,-1.5) {$\tilde O_2$};
    \draw (A1) -- (beta1) -- (B1) -- (O1) -- (A1);
    \draw (A1) -- (gamma1) -- (C1) -- (tildeO1) -- (B1);  
    \draw[preaction={draw=white, -,line width=6pt}] (A1) -- (beta2) -- (B2) -- (O1);
    \draw[preaction={draw=white, -,line width=6pt}] (A1) -- (gamma2) -- (C2) -- (tildeO2) -- (B2);
     \draw (A1) -- (beta2) -- (B2) -- (O1) -- (A1);
     \draw (A1) -- (gamma2) -- (C2) -- (tildeO2) -- (B2);  
 	\end{tikzpicture}
}
\scalebox{0.8}{
  \begin{tikzpicture}
      \node (A) at (-0.7,0) {\scalebox{0.85}{$A$}};
       \node (B) at (2,0) {\scalebox{0.85}{$B$}};
       \node (C) at (5,0) {\scalebox{0.85}{$C$}};
       \node (beta) at (2,1.8) {\scalebox{0.85}{$\beta$}};
       \node (gamma) at (5,1.8) {\scalebox{0.85}{$\gamma$}};
    \node (O) at (1,-1.5) {$O$};
    \node (tildeO) at (4,-1.5) {$\tilde O$};
    \draw (A) -- (beta) -- (B) -- (O) -- (A);
    \draw (A) -- (gamma) -- (C) -- (tildeO) -- (B);  
 		 \end{tikzpicture}
   }
\[ \varphi: \Sigma \to \Delta \]
\caption{
    \label{fig:CounterexampleIDread} 
    On the left, a poset $\Sigma$ and a map $m_1 : \calV \to \ZZ_{\ge 1}$ 
    with $\calV = \Sigma \setminus \aset{O_1, \tilde O_1, \tilde O_2}$;
    right, a poset $\Delta$;
    together, a \combinatorial{} morphism of posets $\varphi$ such that $m_1$ can be extended to $O_1$ but not to $\tilde O_1$ nor $\tilde O_2$.
}
\end{figure}
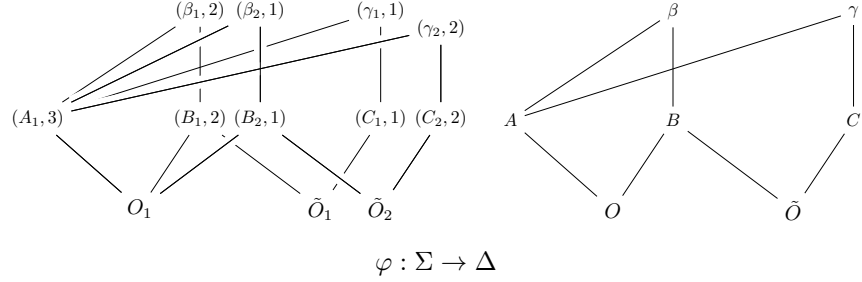

  \begin{ex} 
   \label{ex:CannotAlwaysExtend} 
   It is straightforward to construct balanced maps $m_\calV$ that cannot be extended when $\calV$ is disconnected;
   e.g.~consider the poset $\alpha \leftarrow A \leftarrow O \rightarrow B \rightarrow \beta$,
   the \combinatorial{} morphism $\varphi = \operatorname{id}$,
   the set $\calV = \upset \aset{A, B}$,
   the balanced map $m_\calV(\alpha) = m_\calV(A) = 2$ and $m_\calV(\beta) = m_\calV(B) = 1$.
   It is not possible to extend $m_\calV$ to $O$.

   Moreover, we give an example where $\calV$ is connected, yet extension is not possible. 
   Figure~\ref{fig:CounterexampleIDread} shows a \combinatorial{} morphism  
   $\varphi : \Sigma \to \Delta$ given by $\varphi(A_i) = A$, $\varphi(\beta_i) = \beta$, etc.~and
   a map $m_\calV : \calV \to \Delta$ with $\calV = \Sigma \setminus \aset{O_1, \tilde O_1, \tilde O_2}$ 
   and whose values are given by the second numbers of the pairs in the diagram.
   If we consider $\tilde O_1$, we have  $\varphi(\tilde O_1) = \tilde O$,
   and both $B$ and $C$ cover $\tilde O$.
   Note that
   \[ \sum_{\substack{\eta \in \inv{\varphi}(B) \\ \tilde O_1 \lessdot \eta  }  } m_\calV(\eta) 
   = 2 \ne 1 = 
 \sum_{\substack{\eta \in \inv{\varphi}(C) \\ \tilde O_1 \lessdot \eta  }  } m_\calV(\eta). \]
   Thus, there is no possible value for $\tilde O_1$ to fulfill the balancing condition.
   Same with $\tilde O_2$. On the other hand, $O_1$ can be given the value~3, and this satisfies the balancing condition. 
 \end{ex} 

The following result sheds light on why in Example~\ref{ex:CannotAlwaysExtend} it was possible to extend $m_\calV$ to $O_1$ but not to $\tilde O_1$ nor $\tilde O_2$.


\begin{prop} 
   \label{prop:ExtendIndexMap}
  Let $\varphi : \Sigma \to \Delta$ be a \combinatorial{} morphism 
  and $\calV  \subset  \calW  \subset  \Sigma$ open sets.
  If $\calU_\alpha = (\upset_\Delta \varphi(\alpha) ) \setminus \aset{ \varphi(\alpha) }$ is connected 
  and $\inv{\varphi}(\calU_\alpha) \subset \calV$ for all $\alpha$ in $\calW \setminus \calV$,
  then any balanced map $m_\calV : \calV \to \ZZ_{\ge 1}$ extends to a balanced map $m_{\calW }: \calW \to \ZZ_{\ge 1}$ 
  by setting 
  \begin{align} 
    \label{eq:ExtendingBalancedMap}
    m_{\calW}(\alpha)  = \sum_{\substack{\gamma \in \inv \varphi(\beta) \\ \alpha \lessdot \gamma }} m_\calV(\gamma) 
  \end{align}
  for $\alpha \in \calW \setminus \calV$ and $\beta$ covering $\varphi(\alpha)$. 
  The value is independent of the choice of $\beta$.
\end{prop}

\begin{proof} 
   We show that the value in Equation~\eqref{eq:ExtendingBalancedMap} is independent of the choice of $\beta$;
   this also shows the balancing condition.
   Let $\beta_1, \beta_2$ be elements that cover $\varphi(\alpha)$,
   so they are in~$(\upset_\Delta \varphi(\alpha) ) \setminus \aset{ \varphi(\alpha) } = \calU_\alpha$. 
   Since $\calU_\alpha$ is connected, there is a path from $\beta_1$ to $\beta_2$, 
   which can be completed to a path  
   $ P = \langle \mu_0, \mu_1,  \dots, \mu_{k-1}, \mu_k \rangle \subset  (\upset_\Delta \varphi(\alpha) ) \setminus \aset{ \varphi(\alpha) } $
   such that $\mu_0 = \beta_1$, $\mu_k = \beta_2$  and 
   either  $\mu_{i-1} \lessdot \mu_{i}$ or $\mu_{i-1} \gtrdot \mu_{i}$.
   Since $\inv{\varphi}(\calU_\alpha) \subset \calV$, we have that~$\inv{\varphi}(\mu_i) \subset \calV$.
   Moreover, $\mu_{i-1}$ and $\mu_{i}$ are in $\upset_\Delta \varphi(\alpha)$, 
   so we can apply Lemma~\ref{lm:BalancingFormulas} to consecutive elements of $P$ to obtain the result.
\end{proof}
 
\subsection{Strongly connected posets}
    \label{sub:StronglyConnectedPosets}
We simplify Theorem~\ref{mytheorem:ExtendIndexMap} into a statement that is easier to apply to tropical moduli spaces. 
Recall that a poset $\Sigma$ is \emph{graded} if it admits a rank function, i.e.~a poset morphism $\rk : \Sigma \to (\NN, \le)$ such that $\alpha \lessdot \gamma$ if and only if $\rk \gamma = \rk \alpha + 1$ for all $\alpha, \gamma \in \Sigma$.
We assume that the rank of minimal elements is 0.
For polyhedral complexes $\dim \sigma$ is a natural rank function.
We denote by~$\Sigma(k)$ all the rank-$k$ elements of~$\Sigma$.
Now we recall three standard definitions, and introduce one by us.

\begin{de} 
    \label{de:StrongConnectivity}  
    Let $\Sigma$ be a graded finite poset with rank $\rk$.
   \begin{itemize} 
     \item The dimension of $\Sigma$ is the maximal value of $\rk$.
     \item $\Sigma$ is of \emph{pure dimension d} if $\rk$ is constant on $\max \Sigma$ and the value is~$d$.
     \item $\Sigma$ is \emph{connected in codimension-$k$} if the up-set $\upset_\Sigma \Sigma(d-k)$ is connected, for some $k$ in $\aset{0,\dots,d}$.
     \item $\Sigma$ is \emph{strongly connected} if $\Sigma$ is connected, and for every $\alpha \in \downset_\Sigma \Sigma(d-2)$ the up-set $(\upset_\Sigma \alpha) \setminus \aset{\alpha}$ is connected.
   \end{itemize}
\end{de}

Posets associated to tropical objects typically display pure dimension and being connected in codimension-1.
We show that being strongly connected implies both.


\begin{lm} 
    \label{lm:StronglyConnectedImpliesCodimension1Connected}
   Let $\Sigma$ be a graded poset of dimension $d$. 
   If $\Sigma$ is strongly connected, then $\Sigma$ is connected in codimension-$k$ for any $k$ in $\aset{1,\dots,d}$,
   and its dimension is pure.
\end{lm}

\begin{proof}
  Let $\alpha$ and $\beta$ be elements in $\upset \Sigma(d-k)$.
  Since $\Sigma$ is strongly connected, we have that $\Sigma$ is connected.
  So there is a sequence $P_0 = \langle \gamma_0, \dots, \gamma_q \rangle$ whose consecutive elements are comparable and with $\gamma_0 = \alpha$ and $\gamma_q = \beta$.
  Now, let $q = \min_{\gamma \in P_0} \rk \gamma$.
  If $q \ge d-k$, then $\alpha$ and $\beta$ are already connected in codimension-$k$.
  So assume that $q<d-k$, and let $\gamma_i \in P_0$ be an element with $\rk \gamma_i = q$.
  We have that either $\gamma_{i-1} \ge \gamma_i$ or $\gamma_{i-1} \le \gamma_i$, 
  but the latter would imply $\rk \gamma_{i-1} = q-1$, so it is excluded. 
  The same holds for $\gamma_{i+1}$, so $\gamma_{i-1}$ and $\gamma_{i+1}$ are in $\upset_\Sigma \gamma_i \setminus \aset{\gamma_i}$,
  and since $\rk \gamma_i = q < d-k \le d-1$ and $\Sigma$ is strongly connected, 
  we can connect $\gamma_{i-1}$ with $\gamma_{i+1}$ via a path $Q$ whose minimal rank is $q+1$.
  Thus, we replace $\gamma_i$ in $P_0$ with the path $Q$ to obtain $P_1$.
  Iterating this procedure eliminates all elements in $P_0$ with rank $q$, rank $q+1$, and so on until all have rank at least $d-k$. 
  
  Next, consider $\alpha \in \Sigma$ with $\rk \alpha = q < d$, and $\beta \in \Sigma$ with $\rk \beta= d$.
  As $d-q \ge 1$, the set  $\upset_\Sigma(q)$ is connected,
  so there is a path $P = \langle \alpha, \gamma, \dots, \beta \rangle \subset \upset \Sigma(q)$.
  Again, $\gamma < \alpha$ is excluded by rank reasons, which means $\gamma > \alpha$, so $\alpha$ is not a maximal element of $\Sigma$.
  Hence, all maximal elements of $\Sigma$ have rank~$d$, as desired.
\end{proof}

\begin{re} 
  The poset $\Sigma$ shown in Figure~\ref{fig:CounterexampleIDread} is connected in codimension-1, but is not strongly connected.
\end{re}

So we obtain a situation where the technical conditions of Proposition~\ref{prop:ExtendIndexMap} are met.

\begin{lm} 
    \label{lm:CodimensionkExtendIndexMap} 
  Let $\Sigma$ and $\Delta$ be graded posets,
  $\dtmor : \Sigma \to \Delta$ a \combinatorial{} morphism,
  $\calV = \upset_\Sigma \Sigma(k+1)$ and $\calW = \upset_\Sigma \Sigma(k)$, where $k$ is in $\aset{0, \dots, \rk \Sigma - 2}$.
  For all $\alpha \in \calW \setminus \calV$ we have that $\inv{\varphi}( (\upset_\Delta \varphi(\alpha) ) \setminus \aset{ \varphi(\alpha) } ) \subset \calV$ 
\end{lm}

\begin{proof} 
  Note that $\rk \alpha = k$, because $\alpha$ is in $\calW \setminus \calV$, and that $\dtmor$ preserves rank because $\dtmor$ is \combinatorial{}.
  Hence, $\rk \dtmor(\alpha) = k$ and
  all elements of $(\upset_\Delta \varphi(\alpha) ) \setminus \aset{ \varphi(\alpha) }$ have rank at least $k+1$, 
  and so do all elements of $\inv{\varphi}( (\upset_\Delta \varphi(\alpha) ) \setminus \aset{ \varphi(\alpha) } )$.
  Therefore, the latter set is in $\calV = \upset_\Sigma \Sigma(k+1)$, as desired. 
\end{proof}

\begin{proof}[Proof of Proposition~\ref{myprop:ExtendingMap}]
  By induction, suppose $m_\calV$ is a balanced map defined on $\calV = \upset \Sigma(k+1)$.
  Set $\calW = \upset \Sigma(k)$.
  Since $\Delta$ is strongly connected we have that $\calU_\alpha = (\upset_\Delta \varphi(\alpha) ) \setminus \aset{ \varphi(\alpha) }$ is connected.
  By Lemma~\ref{lm:CodimensionkExtendIndexMap}, we have that $\inv{\varphi}(\calU_\alpha) \subset \calV$.
  Thus, by Theorem~\ref{mytheorem:ExtendIndexMap} the map $m_\calV$ extends to a balanced map $m_\calW$. 
\end{proof}

\subsection{Further connectivity of posets}
    \label{sub:FurtherConnectivityOfPosets}
When certain corestriction of $\dtmor$ is \combinatorial, we can go beyond upward paths in our study of connectivity, and lift general paths.
This allows, in the presence of a connected fibre, to lift the connectivity of $\Delta$ to $\Sigma$.

\begin{lm} 
    \label{lm:LiftingPaths}  
      Let $\varphi : \Sigma \to \Delta$ be a morphism, 
      $\calV \subset \Sigma$ an up-set, 
      $m_\calV$ a balanced map,
      $\alpha$ in $\calV$ an element, and
      $P = \langle \beta_0, \beta_1, \dots, \beta_k \rangle \subset  \varphi(\calV)$ a path with $\beta_0 = \varphi(\alpha)$.
      If $\psi = (\varphi|_{\calV})|^{\varphi(\calV)}$ is \combinatorial{},
      there is a lift $\tilde P = \langle \gamma_0, \gamma_1, \dots, \gamma_q \rangle$ in $\calV$ with $\gamma_0 = \alpha$.
\end{lm}

\begin{proof} 
      Proceeding by induction, suppose that $\beta_i$ lifts to $\gamma_i$ in $\calV$,
      an find a lift $\gamma_{i+1} \in \calV$ of $\beta_{i+1}$
      that is comparable to $\gamma_i$.
      The case $\beta_i \le \beta_{i+1}$ is handled by Lemma~\ref{lm:LiftingUpwardsPaths}. 
      Assume that $\beta_i \ge \beta_{i+1}$. 
      Since $\psi$ is \combinatorial{}, 
      the down-set $\downset_\calV \gamma_i$ is mapped isomorphically to 
      $\downset_{\varphi(\calV)} \varphi(\gamma_i) = \downset_{\varphi(\calV)} \beta_i$.
      As $\beta_{i+1}$ is in $\downset_\Delta \beta_i$ and in $\varphi(\calV)$, 
      it is in $\downset_{\varphi(\calV)} \beta_i$ as well,
      so there is $\gamma_{i+1}$ in $\downset_\calV \gamma_i \subset \calV$ such that $\varphi(\gamma_{i+1}) = \beta_{i+1}$.
\end{proof}

  \begin{figure} \center
		 \begin{tikzpicture}
		\node (A1) at (-1,0) {$A_1$};
		\node (A2) at (0.3,0) {$A_2$};
		\node (B1) at (1.6,0) {$B_1$};
    \node (B2) at (3,0) {$(B_2,1)$};
    \node (C1) at (5,0) {$(C_1,2)$};
    \node (O1) at (2,-1) {$O_1$};
    \node (t1) at (4.5,1.2) {$(\beta_2,1)$};
    \node (t2) at (2.5,1.2) {$(\beta_1,1)$};
		\node (s1) at (-0.6,1.2) {$\alpha_1$};
    \node (s2) at (0.9,1.2) {$(\alpha_2,1)$};
    \draw (B2) -- (O1) -- (A1) -- (s1) -- (B1);
    \draw (B2) -- (t1) -- (C1);  
		\draw (B1) -- (O1) -- (C1) -- (t2); 
    \draw (A2) -- (O1);
		 \draw[preaction={draw=white, -,line width=6pt}] (B1) -- (t2);
     \draw[preaction={draw=white, -,line width=6pt}] (t1) -- (B2) -- (s2) -- (A2);
		 \node (A) at (6,0) {$A$};
		 \node (B) at (8,0) {$B$};
		 \node (C) at (10,0) {$C$};
		 \node (t) at (9,1.2) {$\beta$};
		 \node (s) at (7,1.2) {$\alpha$};
     \node (O) at (8,-1) {$O$};
     \draw (O) -- (A) -- (s) -- (B) -- (t) -- (C) -- (O);
     \draw (O) -- (B);
		 \node (t2) at (4,-1.5) {\scalebox{1.25}{$\Sigma$}};
		 \node (t2) at (5,-1.5) {\scalebox{1.4}{$\rightarrow$}};
		 \node (t) at (6,-1.5) {\scalebox{1.25}{$\Delta$}};
		 \end{tikzpicture}
\caption{
      \label{fig:BalancedYetNotCover} 
      A \combinatorial{} morphism $\varphi : \Sigma \to \Delta$, 
      and a balanced map $m_\calV : \calV \to \ZZ_{\ge 1}$ with $\calV = \upset_\Sigma \aset{B_2, C_1}$ 
      such that $\varphi$ restricted to $\calV$ is not an indexed branched cover.
      } 
\end{figure}

\begin{ex}
      \label{ex:PathThatDoesNotLift}
  We wonder if in Lemma~\ref{lm:LiftingPaths}
  one can omit the condition of having a \combinatorial{} $\psi$.
  It could, for example, follow from a weaker condition such as requiring $\varphi$ to be \combinatorial{}.
  This is not the case, e.g.~in Figure~\ref{fig:BalancedYetNotCover} 
  let $\calV = \upset_\Sigma \aset{B_2, C_1} = \aset{B_2, C_1, \beta_1, \beta_2, \alpha_2}$.
  We have that $\beta_1 \in \calV$, and $\downset_\calV \beta_1 = \aset{\beta_1, C_1}$, 
  but $\downset_{\varphi{\calV}} \psi(\beta_1) = \aset{\beta, C, B}$, 
  so $\psi$ is not \combinatorial{}.
  Moreover, this example has a path that cannot be lifted:
  consider the path $\langle \beta, B \rangle \subset \varphi(\calV)$, there is no lift in $\calV$ starting with $\beta_1$.

  Figure~\ref{fig:BalancedYetNotCover} also underscores the fact that balanced maps restrict to balanced maps for any up-set $\calV$, but the issue is more subtle with indexed branched covers.
  Setting $m_\calV(C_1) = 2$ and $m_\calV$-value 1 for the remaining elements in $\calV$, 
  we get a balanced map $m_\calV$, yet $(\psi, m_\calV)$ is not an indexed branched cover,
  since the count for the fibre over $C$ is 2, 
  and for the fibre over $B$ is~1.
\end{ex}


We now use the results on lifting paths to give a criterion for the connectedness of $\calV$. 

\begin{prop} 
    \label{prop:OneFibreConnection}  
  Let $\varphi : \Sigma \to \Delta$ be a morphism, 
  $\calV \subset \Sigma$ an up-set such that $\psi = (\varphi|_{\calV})|^{\varphi(\calV)}$ is \combinatorial{},
  and $m_\calV : \calV \to \ZZ_{\ge 1}$ a balanced map.
  If $\varphi(\calV)$ is connected and there is $\beta$ in $\varphi(\calV)$ such that $\inv{\varphi}(\beta)$ 
  is connected in $\calV$,
  then $\calV$ is connected.
\end{prop}

\begin{proof} 
  Consider $\alpha \in \calV$.  
  There is a path $P = \langle \varphi(\alpha), \beta_1, \dots, \beta \rangle$ connecting $\varphi(\alpha)$ and~$\beta$,
  since both are in $\varphi(\calV)$.
  By Lemma~\ref{lm:ConnectedImpliesSequence}, there is a lift $\tilde P = \langle \alpha, \nu_1, \dots, \nu_{k-1}, \nu_k \rangle$
  connecting $\alpha$ with fibre $\inv{\varphi}(\beta)$.
  Since the fibre $\inv{\varphi}(\beta)$ is connected, we are done. 
\end{proof}

\subsection{Restricting \combinatorial{} maps}
    \label{sub:RestrictingCombiMaps}
  It remains to prove Proposition~\ref{myprop:LiftingConnectivity}.
  We do this by studying which up-sets $\calV$ satisfy that the restriction $\psi = (\varphi|_{\calV})|^{\varphi(\calV)}$ of any \combinatorial{} $\varphi$ is \combinatorial{} as well.
  Contrast the following result with Example~\ref{ex:PathThatDoesNotLift}.

\begin{lm} 
    \label{lm:RestrictToConnectedComponentFibre}
  Let $\dtmor : \Sigma \to \Delta$ be a \combinatorial{} morphism and $\calU \subset \Delta$ an up-set.
  If $\calV$ is a union of connected components of $\inv \dtmor(\calU)$,
  then  $\psi = (\varphi|_{\calV})|^{\varphi(\calV)}$ is \combinatorial{}.
\end{lm}

\begin{proof} 
  Since $\varphi$ is \combinatorial{}, for every $\alpha \in \calV$ the set $\downset_{\calV} \alpha$ is mapped isomorphically by $\psi$ to $\psi( \downset_{\calV} \alpha)$, 
  and we must show this coincides with $\downset_{\dtmor(\calV)} \psi(\alpha)$.
  Note that $\psi(\downset_{\calV} \alpha)$ equals $\psi(\downset_\Sigma \alpha \cap \calV)$, 
  and $\downset_{\varphi(\calV)} \psi(\alpha)$ equals $\downset_{\Delta} \psi(\alpha) \cap \varphi(\calV)$.
  So we only need to show that  if $\beta \in \downset_{\Delta} \psi(\alpha) \cap \varphi(\calV)$, then $\beta$ is in $\psi(\downset_\Sigma \alpha \cap \calV)$; the other containment is clear.

  As $\dtmor$ is \combinatorial{} and $\beta$ is in $\downset_{\Delta} \dtmor(\alpha)$, there is a unique lift $\gamma$ of $\beta$ such that $\gamma \preceq \alpha$.
  Consider the upwards path $\tilde P = \langle \gamma, \alpha \rangle$, mapped by $\dtmor$ to $P = \langle \dtmor(\gamma), \dtmor(\alpha) \rangle = \langle \beta, \psi(\alpha) \rangle$.
  The path $P$ is in $\dtmor(\calV)$, a subset of $\calU$, because $\beta$ is in $\dtmor(\calV)$.
  Thus, $\tilde P$ is a path in $\inv \dtmor(\calU)$.
  Since $\calV$ is a union of connected comoponents of $\inv \dtmor(\calU)$, any path $\tilde P$ contained in $\inv \dtmor(\calU)$ such that $\tilde P \cap \calV \ne \varnothing$, satisfies that $\tilde P \subset \calV$.
  Hence, $\alpha \in \tilde P \cap \calV$ implies that $\tilde P \subset \calV$, thus $\gamma$ is in $\calV$.
  Therefore, $\beta = \psi(\gamma)$ is in $\psi(\downset_\Sigma \alpha \cap \calV) $, as desired. 
\end{proof} 

The previous result, combined with the fact that a \combinatorial{} morphism preserves the rank function of a graded poset, gives our last proof.

\begin{proof}[Proof of Theorem~\ref{myprop:LiftingConnectivity} ]
  Recall that we are assuming all minimal elements have rank equal to 0.
  Since $\dtmor$ is \combinatorial{}, the set $\downset_\Sigma \alpha$ is mapped isomorphically to $\downset_\Delta \varphi(\alpha)$, which means $\rk_\Sigma \alpha = \rk_\Delta \varphi(\alpha)$. 
  Thus, if we set $\calV = \upset \Sigma(k)$ and $\calU = \upset \Delta(k)$, we have that $\calV = \inv \varphi(\calU)$.
  By Lemma~\ref{lm:RestrictToConnectedComponentFibre}, this gives that  $\psi = (\varphi|_{\calV})|^{\varphi(\calV)}$ is \combinatorial{}.
  So by Proposition~\ref{prop:OneFibreConnection}, if there is $\beta \in \calU$ such that $\inv \varphi(\beta)$ is connected in $\calV$, then $\calV$ is connected and we are done.
\end{proof}

\printbibliography

\end{document}